\documentclass[12pt, reqno]{amsart}
\usepackage{amssymb,latexsym,amsmath,amscd,amsthm,graphicx, color}
\usepackage[all]{xy}
\usepackage{pgf,tikz}
\usepackage{mathrsfs}
\usetikzlibrary{arrows}
\usepackage{pgfplots}
\usepackage[left=0.6 in, top=0.6 in, right=0.6 in, bottom=0.3 in]{geometry}
\usepackage{changepage}
\allowdisplaybreaks
\makeatletter
\newcommand{\setword}[2]{%
	\phantomsection
	#1\def\@currentlabel{\unexpanded{#1}}\label{#2}%
}
\makeatother

\usepackage[linkcolor=blue, urlcolor=blue, citecolor=blue,
colorlinks, bookmarks]{hyperref}

\definecolor{uuuuuu}{rgb}{0.26666666666666666,0.26666666666666666,0.26666666666666666}
\definecolor{xdxdff}{rgb}{0.49019607843137253,0.49019607843137253,1.}
\definecolor{ffqqqq}{rgb}{1.,0.,0.}
\definecolor{ffqqqq}{rgb}{1.,0.,0.}
\definecolor{ffxfqq}{rgb}{1.,0.4980392156862745,0.}

\definecolor{uuuuuu}{rgb}{0.26666666666666666,0.26666666666666666,0.26666666666666666}
\definecolor{qqwuqq}{rgb}{0.,0.39215686274509803,0.}
\definecolor{zzttqq}{rgb}{0.6,0.2,0.}
\definecolor{xdxdff}{rgb}{0.49019607843137253,0.49019607843137253,1.}
\definecolor{qqqqff}{rgb}{0.,0.,1.}
\definecolor{cqcqcq}{rgb}{0.7529411764705882,0.7529411764705882,0.7529411764705882}
\definecolor{sqsqsq}{rgb}{0.12549019607843137,0.12549019607843137,0.12549019607843137}
\definecolor{uuuuuu}{rgb}{0.26666666666666666,0.26666666666666666,0.26666666666666666}
\definecolor{ffqqqq}{rgb}{1,0,0}
\definecolor{xdxdff}{rgb}{0.49019607843137253,0.49019607843137253,1}

\definecolor{yqqqyq}{rgb}{0.5019607843137255,0,0.5019607843137255}
\definecolor{qqqqff}{rgb}{0,0,1}
\definecolor{ffqqqq}{rgb}{1,0,0}
\definecolor{ffqqff}{rgb}{1,0,1}

\theoremstyle{plain}

\newtheorem{theorem}[subsection]{Theorem}
\newtheorem{theo1}[subsubsection]{Theorem}

\newtheorem{prop}[subsection]{Proposition}
\newtheorem{prop1}[subsubsection]{Proposition}
\newtheorem{lemma1}[subsubsection]{Lemma}

\theoremstyle{definition}

\newtheorem{remark}[subsection]{Remark}
\newtheorem{remark1}[subsubsection]{Remark}

\newtheorem{open}[subsubsection]{Open}

\newtheorem{notation}[subsubsection]{Notation}
\newtheorem{note}[subsubsection]{Note}
\newtheorem{note1}[subsection]{Note}

\newcommand{\uu}{\cup}
\newcommand{\ii}{\cap}
\newcommand{\UU}{\bigcup}

\newcommand{\set}[1]{\{#1\}}

\newcommand{\ga}{\beta}
\newcommand{\gb}{\gamma}

\newcommand{\gk}{\kappa}

\newcommand{\gq}{\theta}
\newcommand{\rh}{\sigma}

\newcommand{\tbf}{\textbf}
\newcommand{\tit}{\textit}

\newcommand{\D}[1]{\mathbb{#1}}

\newcommand{\te}{\text}

\newcommand{\im}{\implies}

\begin{document}

\title{Optimal Sets and Quantization Errors under Geometric Constraints for Discrete Distributions}

\address{$^{1, 3}$Department of Mathematics\\
Dr B R Ambedkar National Institute of Technology Jalandhar\\
Jalandhar Punjab, India 144011.}
\address{$^{4}$Department of Mathematics\\
 Indian Institute of Technology Delhi\\
  New Delhi, India 110016.}
\address{$^{2, 5}$School of Mathematical and Statistical Sciences\\
		University of Texas Rio Grande Valley\\
		1201 West University Drive\\
		Edinburg, TX 78539-2999, USA.}

\email{$^1$prabhat65tamrakar@gmail.com, $^3$sbhkatiyar@gmail.com,  $^4$maz248182@iitd.ac.in}
\email{\{$^2$bismark.bimpong01, $^5$mrinal.roychowdhury\}@utrgv.edu}

 \author{$^1$Prabhat Tamrakar}
\author{$^2$Bismark Bimpong}
\author{$^3$S. K. Katiyar}
 \author{$^4$Sayandip Pandit}
  \author{$^5$Mrinal Kanti Roychowdhury}



\date{}
\maketitle

\pagestyle{myheadings}\markboth{Tamrakar, Bimpong, Katiyar, Pandit, Roychowdhury}
{Constrained quantization for discrete distributions}

\begin{abstract}

This paper presents a detailed study of constrained quantization for both finite and infinite discrete probability distributions supported on subsets of the real line. Under specific geometric constraints - namely, a semicircular arc and the union of two sides of an equilateral triangle - we compute constrained optimal sets of $n$-points and the corresponding $n$th constrained quantization errors. For finite discrete distributions, we consider both uniform and nonuniform cases with support on $\{-3, -2, \dots, 3\}$. For infinite discrete distributions, two cases are analyzed: one supported on $\left\{ \frac{1}{n} : n \in \mathbb{N} \right\}$ and the other on the set of natural numbers $\mathbb{N}$. Explicit constructions and numerical computations of optimal quantizers and errors are provided. Furthermore, for the infinite discrete distribution supported on $\mathbb{N}$, we develop a general framework for constrained quantization under the linear constraint $y = mx + c$ and prove that the constrained quantization dimension in this setting is zero. Our results highlight how geometric constraints influence the structure and existence of optimal quantizers and pave the way for further investigations into constrained quantization theory.
\end{abstract}

\section{Introduction}
Quantization is the process of reducing a digital signal with $n $ levels to one with $k $ levels, where $k < n $, or converting a continuous analog signal into a digital signal with $k $ discrete levels. This process is fundamental for representing, processing, storing, and transmitting analog information in digital systems, as well as for achieving data compression. It continues to be a key and actively researched area in information theory and source coding, with numerous applications in technology and engineering, especially in data compression and clustering. For foundational work, see  \cite{DFG}-\cite{Z2} .
More recently, the scope of quantization theory has expanded with the development of constrained and conditional quantization (see \cite{PR1, PR2}). This has led to a classification of the theory into two principal areas: constrained and unconstrained quantization. For a concise overview and practical entry point into both branches, see \cite{PR1}. For a deeper exploration of unconstrained quantization, refer to \cite{GL2}. The present paper focuses specifically on constrained quantization.

Let $ \mu $ be a Borel probability measure on $ \mathbb{R}^2 $ equipped with the Euclidean metric $ d $ induced by a norm $ \|\cdot\| $. 
It is assumed that $ R $ is a nonempty closed subset of $\mathbb{R}^2 $ and that $\ga \subseteq R $ is a locally finite subset of $\mathbb{R}^2 $. This means that the intersection of $\ga $ with any bounded subset of $\mathbb{R}^2 $ is finite. This suggests that $ \ga $ is closed and countable.
With regard to the set $ \ga \subseteq R$, the distortion error for the measure $ \mu $ is defined by
\[
Q(\mu; \ga ) = \int \min_{b \in \ga } \|x - b\|^2 \, d\mu(x).
\]
Then, for $ n \in \mathbb{N} $, the \textit{$ n$th constrained quantization error} for $ \mu $ with respect to the set $ \mathit{R} $ is defined by
\begin{equation} \label{Vr}
Q_n := Q_n(\mu) = \inf \left\{ Q(\mu; \ga ) : \ga \subseteq R,\ 1 \leq \operatorname{card}(\ga ) \leq n \right\},
\end{equation}
where $ \operatorname{card}(A) $ denotes the cardinality of a set $ A $. Let $Q_{n}(\mu)$ be a strictly decreasing sequence, and write $Q_{\infty}(\mu):=\mathop{\lim}\limits_{n\to \infty} Q_{n}(\mu)$.
Then, the number $D(\mu)$ defined by
\begin{align*}
D(\mu):=\mathop{\lim}\limits_{n\to \infty}  \frac{2\log n}{-\log (Q_{n}(\mu)-Q_{\infty}(\mu))}
\end{align*}
if it exists,
 is called the \tit{constrained quantization dimension} of $\mu$ and is denoted by $D(\mu)$. The pace at which the given measure of the constrained quantisation error converges as $n$ goes to infinity is measured by the constrained quantisation dimension.
For any $\gk>0$, the two numbers
								\[\liminf_{n\to \infty} n^{\frac 2 \gk}  (Q_{n}(\mu)-Q_{\infty}(\mu)) \te{ and } \limsup_{n\to \infty}  n^{\frac 2 \gk}(Q_{n}(\mu)-Q_{\infty}(\mu))\] are, respectively, called the \tit{$\gk$-dimensional lower} and the \tit{upper} constrained quantization coefficients for $\mu$. If the $\gk$-dimensional lower and the upper constrained quantization coefficients for $\mu$ are finite and positive, then the constrained quantization dimension $D(\mu)$ of $\mu$ exists and $\gk$ equals $D(\mu)$ (see \cite{PR1}).

Let $ U $ be the largest open subset of $ \mathbb{R}^2 $ such that $ \mu(U) = 0 $. Then, the set $ \mathbb{R}^2 \setminus U $ is called the \textit{support of $ \mu $}, denoted by $ \operatorname{supp}(\mu) $.
For a locally finite set $ \ga \subseteq \mathbb{R}^2 $ and a point $ b \in \ga $, define
\[
M(b | \ga ) = \left\{ x \in \mathbb{R}^2 : d(x, b) = \min_{d \in \ga } d(x, d) \right\}.
\]
The set $ M(b | \ga ) $ is called the \textit{Voronoi region} generated by $ b \in \ga $. The collection
\[
\left\{ M(b | \ga ) : b \in \ga \right\}
\]
is called the \textit{Voronoi diagram} (or \textit{Voronoi tessellation}) of $ \mathbb{R}^2 $ with respect to $ \ga $.
A set $ \ga $ for which the infimum in \eqref{Vr} is attained and $ \operatorname{card}(\ga ) \leq n $ is called a \textit{constrained optimal set of $ n $-points} for $ \mu $ with respect to the constraint $ R $.
 If $ R = \mathbb{R}^2 $, i.e., there is no constraint, then the notion of constrained quantization error reduces to the classical definition of unconstrained quantization error, often simply referred to as the quantization error (see \cite{GL2}). In the unconstrained setting, if the support of $ \mu $ contains at least $ n $ elements, then each element in an optimal set of $ n $-points corresponds to the conditional expectation over its own Voronoi region. Due to this property, an optimal set of $ n $-points in unconstrained quantization is also referred to as an optimal set of $ n $-means.

\subsection{Delineation}
In this paper, we study constrained quantization for both finite and infinite discrete distributions. Section~\ref{sec0} provides the necessary preliminaries. We begin by considering two specific geometric constraints: one defined by a semicircular arc with base on the line segment $[-3, 3]$, and the other formed by the union of two sides of an equilateral triangle with the same base. Under these constraints, we determine the constrained optimal sets of $n$-points and the corresponding $n$th constrained quantization errors for all possible $n \in \mathbb{N}$.
Section~\ref{sec1} is devoted to the case of a uniform discrete distribution supported on the set $\{-3, -2, -1, 0, 1, 2, 3\}$, and is divided into two subsections corresponding to the two constraints. Section~\ref{sec2} mirrors this structure but addresses a nonuniform discrete distribution supported on the same set.
In Section~\ref{sec3}, we consider an infinite discrete distribution supported on the set $\left\{\frac{1}{n} : n \in \mathbb{N} \right\}$. The first subsection deals with the semicircular arc with base $[0, 1]$ as the constraint, for which we compute the constrained optimal sets $n$-points and the corresponding $n$th constrained quantization errors for all possible $n \in \mathbb{N}$. In the second subsection, using the constraint formed by two sides of an equilateral triangle with base $[0, 1]$, we perform the same computations for $n$ ranging from $1$ to $2000$. An open problem is also mentioned at the end of this subsection. In Section~\ref{sec4}, we calculate the constrained optimal sets of $n$-points and the $n$th constrained quantization errors for all $n\in \D N$ for an infinite discrete distribution with support the set $\D N$ of natural numbers and the constraint $L$ as a general straight line given by $y=mx+c$. In this section, we also show that the constrained quantization dimension of the measure $\mu$ exists and equals zero. This leads us to investigate some open question as it is mentioned in Remark~\ref{remark467}. Finally in Section~\ref{sec5}, we give the conclusion of the paper and some direction about the future work.

\section{Preliminaries} \label{sec0}
For any two elements $(u_1, u_2)$ and $(v_1, v_2)$ in $\mathbb{R}^2$, the squared Euclidean distance is given by
\[
\rh((u_1, u_2), (v_1, v_2)) := (u_1 - v_1)^2 + (u_2 - v_2)^2.
\]
If the Voronoi regions of two elements $x, y \in \beta$ share a boundary, they are referred to as \emph{adjacent}. Let $z$ be a point on the boundary between the Voronoi regions of $x$ and $y$. Since this boundary is the perpendicular bisector of the line segment connecting $y$ and $y$, it satisfies the \emph{canonical equation}
\[
\rh(x, z) - \rh(y, z) = 0.
\]

Notice that any real number $x \in \mathbb{R}$ can be identified with the element $(x, 0) \in \mathbb{R}^2$. Using this identification, we define a function
\[
\rh : \mathbb{R} \times \mathbb{R}^2 \to [0, \infty) \te{ such that } \rh(x, (u, v)) = (x - u)^2 + v^2,
\]
which represents the squared Euclidean distance between $x \in \mathbb{R}$ and $(u, v) \in \mathbb{R}^2$.
 \begin{remark} \label{rem1}
 Let $ \ga_n$ be an optimal set of $ n $-points for a Borel probability measure $ \mu $. If a point in the support of $ \mu $ lies on the common boundary of two Voronoi regions generated by elements of $ \ga_n$, then the point is assigned to one of the Voronoi regions.
 \end{remark}

 \begin{remark} \label{rem2}
For $n\in\D N$, let $\ga_n$ denote an optimal set of $n$-points for a Borel probability measure $\mu$ with respect to a constraint $R$ in $\D R^2$. Then, for $\ga=\{b_1,b_2,\ldots\}\subseteq R$, let us define the sets $A_{b_i|\beta}$ for $a_i\in \ga$ as follows:
	\begin{equation}\label{Megha1000}
		\begin{aligned}
			A_{b_{i}|\beta}=\left\{\begin{array} {ll}
				M(b_1|\ga) & \te{ if } i=1,\\
				M(b_i|\ga)\setminus \UU_{k<i} M(b_k|\ga) & \te{ if } i\geq 2,
			\end{array}\right.
		\end{aligned}
	\end{equation}
where $ M(b_i|\ga ) $ is the \textit{Voronoi region} generated by $ b_i \in \ga $.
	The set $\set{A_{b_{i}|\beta} : b_i\in \ga}$ is called the \tit{Voronoi partition} of $\D R^k$ with respect to the set $\ga$ (and $R$).  Then, for some positive integer $n$, an optimal set $\ga_n$  does not exist, by that it is meant that there is no $\ga\subseteq R$ containing at least $n$ elements such that
	\[\mu(A_{b|\ga})>0 \te{ for each } b \in \ga,\]
and vice versa. 			
 \end{remark}

 \begin{remark} \label{rem3}
In view of Remark~\ref{rem2}, for some positive integer $n$ if a constrained optimal set $\ga_n$ does not exist, then $\ga_{k}$  does not exist for any positive integer $k\geq n$.
 \end{remark}
\begin{note} \label{note1}
Let $\mu$ be a finite discrete distribution with support $\set{r_1, r_2, \cdots, r_N}$ for some positive integer $N\geq 2$ associated with a probability vector $(p_1, p_2, \cdots, p_N)$, i.e., the probability mass function $f$ for $\mu$ is given by
\[f(r_j)=\left\{\begin{array}{cc}
 p_j,  & \te{ where } j\in \set{1, 2, 3, \cdots, N},   \\
 0 & \te{otherwise}.
 \end{array}
\right.\]
Let $a, b\in \D R$ be such that $a\leq r_1<r_2<\cdots<r_N\leq b$. Let us now consider an equilateral triangle $ABC$, whose base $AB$ is the segment $[a, b]$, i.e., the coordinates of $A$ and $B$ are, respectively, given by $(a, 0)$ and $(b, 0)$. Let $R:=R_1\uu R_2$, where $R_1$ and $R_2$, respectively, represent the two sides $AC$ and $BC$. Let the equations of the lines $AC$ and $BC$ are, respectively, given by $y=m_1x+c_1$ and $y=m_2x+c_2$. Let $\ga_n$ be a constrained optimal set of $n$-points for $\mu$ with respect to the constraint $R$ for some positive integer $n$. Let $\te{card}(\ga_n\ii R_1)=\ell$ and $\te{card}((\ga_n\ii R_2)\setminus (\ga_n\ii R_1))=m$. Then, we have $\ell+m=n$. Then, there exists a set of $n$ distinct real number $b_j$ for $1\leq j\leq n$ with $b_1<b_2<\cdots <b_n$ such that
\begin{align*}
\ga_n\ii R_1=\set{(b_j,  m_1 b_j+c_1) : 1\leq j\leq \ell} \te{ and } \ga_n\ii R_2=\set{(b_j,  m_2 b_j+c_2) : 1\leq j\leq m}
\end{align*}
If an optimal set $\ga_n$ exists, then by Remark~\ref{rem2}, the Voronoi region of each element in $\ga_n$ must contain at least an element from the support of $\mu$.
Let
\begin{align*}\te{card}(M((b_j, m_1 b_j+c_1)|\ga_n)) &=n_j \te{ for } 1\leq j\leq \ell, \te{ and } \\
\te{card}(M((b_j, m_2 b_j+c_2)|\ga_n) \setminus  M(&(b_j, m_1 b_j+c_1) |\ga_n)) =n_{\ell+j} \te{ for } 1\leq j\leq  m.
\end{align*}
 Hence, we can assume that with each $\ga_n$ there exists a vector $(n_1, n_2, \cdots, n_{\ell}, n_{\ell+1}, n_{\ell+2}, \cdots, n_{\ell+m})$ associated with $\ga_n$ such that
\[n_1+n_2+\cdots +n_{\ell+m}=N.\]
The vector $(n_1, n_2, \cdots, n_{\ell}, n_{\ell+1}, n_{\ell+2}, \cdots, n_{\ell+m})$ associated with an optimal set of $n$-points will be referred to as a \tit{canonical vector}.
If $Q_n$ is the corresponding $n$th constrained quantization error, then
\begin{equation} \label{eq890}
 \begin{aligned} Q_n& :=Q_{n_1, n_2, \cdots, n_{\ell+m}}(\mu; \ga)=\sum_{j=1}^{n_1}p_j \rh(r_j, (b_1, m_1b_1+c_1))+\sum_{j=n_1+1}^{n_1+n_2} p_j \rh(r_j, (b_2, m_1b_2+c_1))\\
& + \cdots+\sum_{j=n_1+n_2+\cdots+n_{\ell-1}+1}^{n_1+n_2+\cdots+n_{\ell-1}+n_{\ell}}p_j \rh(r_j, (b_{\ell}, m_1b_{\ell}+c_1))\\
&  +\sum_{j=n_1+n_2+\cdots+n_{\ell-1}+n_{\ell}+1}^{n_1+n_2+\cdots+n_{\ell-1}+n_{\ell+1}}p_j \rh(r_j, (b_{\ell+1}, m_2b_{\ell+1}+c_2)) + \cdots\\
&+\sum_{j=n_1+n_2+\cdots+n_{\ell-1}+n_{\ell}+\cdots+n_{\ell+m-1}+1}^{n_1+n_2+\cdots+n_{\ell-1}+n_{\ell+1}+\cdots+n_{\ell+m-1}+n_{\ell+m}}p_j \rh(r_j, (b_{\ell+m}, m_2b_{\ell+m}+c_2)).
\end{aligned}
\end{equation}
\qed
\end{note}

In the following sections we give the main results of the paper.

\section{Constrained quantization for a finite uniform discrete distribution} \label{sec1}
Let $\mu$ be a finite uniform distribution with support $\set{-3, -2, -1, 0, 1, 2, 3}$. Then, the probability mass function $f$ for $\mu$ is given by
 \[f(j)=\left\{\begin{array}{cc}
 \frac 17,  & \te{ where } j\in \set{-3, -2, -1, 0, 1, 2, 3},   \\
 0 & \te{otherwise}.
 \end{array}
 \right.\]
 Comparing with the list $\set{r_1, r_2, \cdots, r_N}$ as defined in Note~\ref{note1}, we have
 \[r_j=-4+j \te{ for } j\in \set{1, 2, 3, 4, 5,  6, 7}.\]
In the following two subsections we calculate the constrained optimal set of $n$-points for the probability distribution $\mu$ for all possible values of $n\in \D N$ with respect to two different constraints.

 \subsection{Constrained quantization when the constraint is a semicircular arc} \label{sub11}
  In this case let us take the constraint $R$ as the upper semicircular arc of the circle $x^2+y^2=9$, i.e.,
 \[R:=\set{(x, y) : x^2+y^2=9 \te{ and } y\geq 0}.\]
 Notice that $R$ can also be represented by
 \[R=\set{(3\cos \gq, 3\sin \gq): 0\leq \gq\leq \pi}.\]
We know that an optimal set of one-point always exists. For any $n\geq 2$, since the boundary of the Voronoi regions of any two optimal elements, in this case, passes through the center of the circle, from the geometry, we see that among $n$ Voronoi regions, only two Voronoi regions contain elements from the support of $\mu$, i.e., only two Voronoi regions have positive probability. Hence, in view of Remark~\ref{rem2} and Remark~\ref{rem3}, the optimal sets of $n$-points exist only for $n=1$ and $n=2$, and they do not exist for any $n\geq 3$.

We now calculate the optimal sets of one-point and two-points in the following propositions.

\begin{prop1}
	Any element on the circle forms an optimal set of one-point with constrained quantization error $Q_1=13$.
\end{prop1}

\begin{proof} Let $(3\cos \gq, 3\sin \gq)$ be an element on the circle. Then, the distortion error for $\mu$ with respect to the element is given by
	\begin{align*}
		Q(\mu; \set{(3\cos \gq, 3\sin \gq)})=&\sum_{j=-3}^3\frac 17\rh(j, (3\cos \gq, 3\sin \gq))=13
		\end{align*}
	which does not depend on $\gq$. Hence, any element on the circle forms an optimal set of one-point with constrained quantization error $Q_1=13$.
\end{proof}

\begin{prop1} \label{prop110}
	The set $\set{(-3, 0), (3, 0)}$ forms an optimal set of two-points with constrained  quantization error $Q_2=\frac {19}{7}$.
\end{prop1}
\begin{proof} Let $\set{(3\cos \gq_1, 3\sin \gq_1), (3\cos \gq_2, 3\sin \gq_2)}$ be an optimal set of two-points, where $0\leq \gq_1<\gq_2\leq \pi$.
Since the boundary of any two elements on the circle passes through the center of the circle, in an optimal set of two-points, one Voronoi region will contain the left half, and the other Voronoi region will contain the right half of the support of $\mu$. Notice that the support element $0$ falls on the boundary of the two Voronoi regions. In that case, as mentioned in Remark~\ref{rem1}, we can assume that  $0$ belongs to the Voronoi region of $(3\cos \gq_2, 3\sin \gq_2)$. Then, the distortion error is given by
	\begin{align*}
		&Q(\mu; \set{(3\cos \gq_1, 3\sin \gq_1), (3\cos \gq_2, 3\sin \gq_2)})\\
&= \sum_{j=-3}^{0}\frac 17\rh(j, (3\cos \gq_2, 3\sin \gq_2))+\sum_{j=1}^{3}\frac 17\rh(j, (3\cos \gq_1, 3\sin \gq_1))\\
&=\frac{1}{7} \left(-36 \cos \gq_1+36 \cos \gq_2+91\right)
		\end{align*}
which is minimum if $\gq_1=0$ and $\gq_2=\pi$, and the minimum value is $\frac {19}{7}$. Thus, the set $\set{(-3, 0), (3, 0)}$ forms an optimal set of two-points with constrained  quantization error $Q_2=\frac {19}{7}$.
\end{proof}

%

\subsection{Constrained quantization when the constraint is the two sides of an equilateral triangle} \label{sub12}  Let us consider the equilateral triangle with vertices $A(-3, 0)$, $B(3, 0)$, and $C(0, 3\sqrt 3)$. Let $R_1$ and $R_2$ represent the sides $AC$ and $BC$, respectively. In this subsection, we take the constraint as $R:=R_1\uu R_2$. Notice that $R_1$ and $R_2$ can be represented as follows:
\begin{align}
R_1=\set{(x, \sqrt{3}(x + 3)) : -3\leq x\leq 0} \te{ and } R_2=\set{(x, \sqrt{3}(-x + 3)) : 0\leq x\leq 3}.
\end{align}

 \begin{prop1}\label{propD1}
An optimal set of one-point is given by $\big\{(-\frac{9}{4},\frac{3 \sqrt{3}}{4})\big\}$ or $\big\{(\frac{9}{4},\frac{3 \sqrt{3}}{4})\big\}$ with constrained quantization error $Q_1=\frac{43}{4}$.
\end{prop1}

\begin{proof} Let $\ga:=\set{(a, b)}$ be an optimal set of one-point. The following two cases can happen:

\tit{Case~1. $(a, b)\in R_1$.}

In this case, $b= \sqrt{3}(a + 3)$ for some $-3\leq a\leq 0$.  Then, the distortion error is given by
\[Q(\mu; \ga)=\sum_{j=-3}^3\frac 17 \rh(j, (a, b))=4 a^2+18 a+31,\]
the minimum value of which is $\frac{43}{4}$, and it occurs when $a=-\frac 94$.

\tit{Case~2. $(a, b)\in R_2$.}

In this case, $b=\sqrt{3}(-a + 3)$ for some $0\leq a\leq 3$. Then, the distortion error is given by
\[Q(\mu; \ga)=\sum_{j=-3}^3\frac 17 \rh(j, (a, b))=4 a^2-18 a+31,\]
the minimum value of which is $\frac{43}{4}$, and it occurs when $a=\frac 94$.

Hence, an optimal set of one-point is given by either $\big\{(-\frac{9}{4},\frac{3 \sqrt{3}}{4})\big\}$ or $\big\{(\frac{9}{4},\frac{3 \sqrt{3}}{4})\big\}$ with constrained quantization error $\frac{43}{4}$, which is the proposition.
\end{proof}

\begin{prop1} \label{propD2}
An optimal set of two-points is $\big\{(-\frac{11}{4}, \frac{\sqrt{3}}{4}), (\frac{21}{8},  \frac{3 \sqrt{3}}{8})\big\}$, or $\big\{(-\frac{21}{8},\frac{3 \sqrt{3}}{8}), (\frac{11}{4},\frac{\sqrt{3}}{4})\big\}$ with constrained quantization error $Q_2=\frac{16}{7}.$
\end{prop1}

\begin{proof}

Let $\ga_2$ be an optimal set of two-points. The following cases can arise:

\tit{Case~1: Both the elements are on $R_1$.}

 Then, we can write
 \[\ga_2:=\set{(a_1, \sqrt{3}(a_1 + 3)), (a_2, \sqrt{3}(a_2 + 3))}, \te{ where } -3\leq a_1<a_2\leq 0.\]
Let $(n_1, n_2)$ be a canonical vector associated with $\ga_2$. Then, all possible values $(n_1, n_2)$ are given by \[\{(1,6),\ (2,5),\ (3,4),\ (4,3),\ (5,2),\ (6,1)\}.\]
Using the formula \eqref{eq890}, we see that the $Q_{n_1, n_2}$ is minimum if $(n_1, n_2)=(3,4) \te{ or } (4,3)$, and the minimum value is given by
\[Q_{3,4}=Q_{4,3}=10.\]

\tit{Case~2: Both the elements are on $R_2$.}

This case is the reflection of Case~1 with respect to the line $x=0$.

\tit{Case~3: One element is on $R_1$ and one element is on $R_2$.}

 In this case, we can write $\ga_2=\set{(a_1, \sqrt{3}(a_1 + 3)), (a_2, \sqrt{3}(-a_2 + 3))}$, where $-3\leq a_1\leq 0 \leq a_2\leq 3$ and $a_1< a_2$.
Then, all possible values of the canonical vector $(n_1, n_2)$ are given by \[\{(1,6),\ (2,5),\ (3,4),\ (4,3),\ (5,2),\ (6,1)\}.\]
Using the formula \eqref{eq890}, we see that the $Q_{n_1, n_2}$ is minimum if $(n_1, n_2)=(3,4) \te{ or } (4,3)$, and the minimum value is given by
\[Q_{3,4}=Q_{4,3}=\frac{16}{7},\]
and it occurs when $(a_1= -\frac{11}{4} \te{ and } a_2=\frac{21}{8})$,  or $(a_1= -\frac{21}{8} \te{ and } a_2=\frac{11}{4})$.

Comparing Case~1 through Case~3, we can deduce that an optimal set of two-points is given by  $\set{(-\frac{11}{4}, \frac{\sqrt{3}}{4}), (\frac{21}{8},  \frac{3 \sqrt{3}}{8})}$, or $\set{(-\frac{21}{8},\frac{3 \sqrt{3}}{8}), (\frac{11}{4},\frac{\sqrt{3}}{4}) }$ with constrained quantization error $Q_2=\frac{16}{7}.$
\end{proof}

 \begin{prop1} \label{propD1}
An optimal set of three-points is given by
$$  \Big\{  \big(-\frac{23}{8},\frac{\sqrt{3}}{8}\big),\big(-\frac{19}{8},\frac{5 \sqrt{3}}{8}\big), \big(\frac{11}{4},\frac{\sqrt{3}}{4}\big) \Big\} $$ with constrained quantization error $Q_3=\frac{15}{7}$.
\end{prop1}

\begin{proof}
Let $\ga_3$ be an optimal set of three-points.
  The following two cases can arise:

\tit{Case~1: All the elements are on one side.}

In this case, without any loss of generality, we can assume that all the elements are on $R_1$. Then, we can write
\[\ga_3=\set{(a_i, \sqrt 3 (a_i+3)) : 1\leq i\leq 3}, \te{ where } -3\leq a_1<a_2<a_3\leq 0.\]
If $(n_1, n_2, n_3)$ be a canonical vector, then $(n_1, n_2, n_3)$ can be an  element  from the set
\begin{align*}
&\set{(1,1,5),\ (1,2,4),\ (1,3,3), (1,4,2),\ (1,5,1),   (2,1,4),\ (2,2,3),\ (2,3,2),\ (2,4,1),  \\
&
 (3,1,3),\ (3,2,2),\ (3,3,1),    (4,1,2),\ (4,2,1),\ (5,1,1)}.
 \end{align*}

 Considering all the possible values of $Q_{n_1, n_2, n_3}(\mu; \ga)$, we see that the minimum value of $Q_{n_1, n_2, n_3}(\mu; \ga)$ is $\frac{69}{7}$ and it occurs when
 \[(n_1, n_2, n_3)=(2,2,3), (2,3,2), \te{ or }(3,2,2).\]

\tit{Case~2: Two elements are on one side and one element on another side.}

Without any loss of generality, we can assume that two elements are on side $R_1$ and one element on side $R_2$. Then, we can write
\[\ga_3=\set{(a_1, \sqrt 3 (a_1+3)),  (a_2, \sqrt 3 (a_2+3)), (a_3, \sqrt 3 (-a_3+3))}, \te{ where } -3\leq a_1<a_2\leq 0\leq a_3\leq 3.\]
 If $(n_1, n_2, n_3)$ be a canonical vector, then $(n_1, n_2, n_3)$ can be an  element  from the set
\begin{align*}
&\set{(1,1,5),\ (1,2,4),\ (1,3,3), (1,4,2),\ (1,5,1),   (2,1,4),\ (2,2,3),\ (2,3,2),\ (2,4,1),  \\
&
 (3,1,3),\ (3,2,2),\ (3,3,1),    (4,1,2),\ (4,2,1),\ (5,1,1)}.
 \end{align*}
  Considering all the possible values of $Q_{n_1, n_2, n_3}(\mu; \ga)$, we see that the minimum value of $Q_{n_1, n_2, n_3}(\mu; \ga)$ is $\frac{15}{7}$ and it occurs when
$(n_1, n_2, n_3)=(2,2,3),$
 and $a_1=-\frac{23}{8}, \, a_2=-\frac{19}{8}$, and $a_3=\frac{11}{4}$.

Taking into account, Case~1 and Case~2, we deduce that the set $\big\{ (-\frac{23}{8},\frac{\sqrt{3}}{8}),(-\frac{19}{8},\frac{5 \sqrt{3}}{8}), (\frac{11}{4},\frac{\sqrt{3}}{4}) \big\}$ forms an optimal set of three-points with constrained quantization error $Q_3=\frac{15}{7}$.
 \end{proof}

 \begin{remark1}
Due to the maximum symmetry of the equilateral triangle with respect to the median passing through the vertex $C$, there is another optimal set of three-points which is the reflection with respect to the median with the same constrained quantization error $Q_3=\frac{15}{7}$.
 \end{remark1}

  \begin{prop1}
The set $\Big\{(-3,0),  (-\frac{21}{8},\frac{3 \sqrt{3}}{8}), (\frac{19}{8},\frac{5 \sqrt{3}}{8}), (\frac{23}{8},\frac{\sqrt{3}}{8})\Big\}$ forms an optimal set of four-points with constrained quantization error $Q_4= \frac{117}{56}$, and there are four different optimal sets of four-points.
\end{prop1}
 \begin{proof}
 Let $\ga_4$ be an optimal set of four-points. Due to the maximum symmetry of the equilateral triangle with respect to the median passing through the vertex $C$, we can assume that $\ga_4$ will contain two elements from $R_1$ and two elements from $R_2$. In fact, it is not difficult to show it.  Then, we can write
\[\ga_4=\set{(a_1, \sqrt 3 (a_1+3)),  (a_2, \sqrt 3 (a_2+3)), (a_3, \sqrt 3 (-a_3+3)), (a_4, \sqrt 3 (-a_4+3))},\]
 where $-3\leq a_1<a_2\leq 0\leq a_3<a_4\leq 3$.
 If $(n_1, n_2, n_3,n_4)$ be a canonical vector, then $(n_1, n_2, n_3, n_4)$ can be an element from the set
\begin{align*}
&\set{(1,1,1,4\},(1,1,2,3),(1,1,3,2),(1,1,4,1),(1,2,1,3),(1,2,2,2),(1,2,3,1),(1,3,1,2),(1,3,2,1),\\
& (1,4,1,1),(2,1,1,3),(2,1,2,2),(2,1,3,1),(2,2,1,2),(2,2,2,1),(2,3,1,1),(3,1,1,2),(3,1,2,1),\\
&(3,2,1,1),(4,1,1,1)}.
 \end{align*}
  Considering all the possible values of $Q_{n_1, n_2, n_3, n_4}(\mu; \ga)$, we see that the minimum value of $Q_{n_1, n_2, n_3,n_4}(\mu; \ga)$ is $\frac{117}{56}$, and it occurs when
 \[(n_1, n_2, n_3, n_4)=(1,2,2,2), (2,1,2,2), (2,2,1,2), \te{ or } (2,2,2,1).\]
 We see that
 \begin{equation} \label{eq90}
\begin{aligned}
 (n_1, n_2, n_3, n_4)&=(1,2,2,2) \im  a_1= -3, a_2= -\frac{21}{8}, a_3= \frac{19}{8}, a_4= \frac{23}{8},\\
  (n_1, n_2, n_3, n_4)&=(2, 1, 2, 2) \im  a_1= -\frac{23}{8}, a_2= -\frac{5}{2}, a_3= \frac{19}{8}, a_4= \frac{23}{8},\\
   (n_1, n_2, n_3, n_4)&=(2, 2, 1, 2) \im  a_1= -\frac{23}{8}, a_2= -\frac{19}{8}, a_3= \frac{5}{2}, a_4= \frac{23}{8},\\
   (n_1, n_2, n_3, n_4)&=(2, 2, 2, 1) \im  a_1= -\frac{23}{8}, a_2= -\frac{19}{8}, a_3= \frac{21}{8},a_4= 3.
\end{aligned}
\end{equation}
Choose $(n_1, n_2, n_3, n_4)=(1,2,2,2)$. Then, the set $\big\{(-3,0),  (-\frac{21}{8},\frac{3 \sqrt{3}}{8}), (\frac{19}{8},\frac{5 \sqrt{3}}{8}), (\frac{23}{8},\frac{\sqrt{3}}{8})\big\}$ forms an optimal set of four-points with constrained quantization error $Q_4= \frac{117}{56}$. By \eqref{eq90}, it follows that there are four different optimal sets of four-points.
 \end{proof}

\begin{prop1}
The set $\Big\{(-3,0),  (-\frac{11}{4},\frac{\sqrt{3}}{4}),  (-\frac{5}{2},\frac{\sqrt{3}}{2}),(\frac{19}{8},\frac{5 \sqrt{3}}{8}), (\frac{23}{8},\frac{\sqrt{3}}{8})\Big\}$ forms an optimal set of five-points with constrained quantization error $Q_5= \frac{29}{14}$, and there are seven different optimal sets of five-points.
\end{prop1}

\begin{proof}
Let $\ga_5$ be an optimal set of five-points. Due to the maximum symmetry of the equilateral triangle with respect to the median passing through the vertex $C$, without any loss of generality, we can assume that $\ga_5$ contains three elements from $R_1$ and two elements from $R_2$. In fact, it is not difficult to show it.  Then, we can write
\[\ga_5=\set{(a_1, \sqrt 3 (a_1+3)),  (a_2, \sqrt 3 (a_2+3)), (a_3, \sqrt 3 (a_3+3)),  (a_4, \sqrt 3 (-a_4+3)), (a_5, \sqrt 3 (-a_5+3)) },\]
 where $-3\leq a_1<a_2<a_3\leq 0\leq a_4<a_5\leq 3$.
Then, the canonical vector $(n_1, n_2, n_3, n_4, n_5)$ can be an element from the following list of vectors
\begin{align*}
\set{& (3,1,1,1,1),(1,3,1,1,1),(1,1,3,1,1),(1,1,1,3,1),(1,1,1,1,3),(2,2,1,1,1),(2,1,2,1,1),(2,1,1,2,1), \\
&(2,1,1,1,2),(1,2,2,1,1),(1,2,1,2,1),(1,2,1,1,2),(1,1,2,2,1),(1,1,2,1,2),(1,1,1,2,2)}.
\end{align*}
   Considering all the possible values of $Q_{n_1, n_2, n_3, n_4, n_5}(\mu; \ga)$, we see that $\frac{29}{14}$ is the minimum value of $Q_{n_1, n_2, n_3, n_4, n_5}(\mu; \ga)$, and it occurs when
\begin{align*} (n_1, n_2, n_3, n_4, n_5)&=(1, 1, 1, 2, 2), (1, 1, 2, 1, 2), (1, 1, 2, 2, 1),(1, 2, 1, 1, 2), (1, 2, 1, 2, 1), (2, 1, 1, 1, 2), \\
& \te{ or } (2, 1, 1, 2, 1),
\end{align*}
i.e., $Q_{n_1, n_2, n_3, n_4, n_5}(\mu; \ga)$ is minimum at seven different canonical vectors. Notice that
\[(n_1, n_2, n_3, n_4, n_5)=(1, 1, 1, 2, 2) \te{ implies that } a_1=-3, \, a_2= -\frac{11}{4}, \, a_3= -\frac{5}{2},\, a_4= \frac{19}{8}, \te{ and } a_5= \frac{23}{8}.\]
Thus, we deduce that the set $\big\{(-3,0),  (-\frac{11}{4},\frac{\sqrt{3}}{4}),  (-\frac{5}{2},\frac{\sqrt{3}}{2}),(\frac{19}{8},\frac{5 \sqrt{3}}{8}), (\frac{23}{8},\frac{\sqrt{3}}{8})\big\}$ forms an optimal set of five-points with constrained quantization error $Q_5= \frac{29}{14}$, and there are seven different optimal sets of five-points.
 \end{proof}

 Let us now state the following two propositions. The proofs are similar to the previous propositions.
  \begin{prop1}
The set $\Big\{(-\frac{23}{8},\frac{\sqrt{3}}{8}), (-\frac{5}{2},\frac{\sqrt{3}}{2}), (-\frac{9}{4},\frac{3 \sqrt{3}}{4}), (\frac{5}{2},\frac{\sqrt{3}}{2}), (\frac{11}{4},\frac{\sqrt{3}}{4}), (3,0)\Big\}$ forms an optimal set of six-points with constrained quantization error $Q_6=\frac{115}{56}$, and there are six different optimal sets of six-points.
\end{prop1}

\begin{prop1}
The set $\Big\{(-3,0), (-\frac{11}{4},\frac{\sqrt{3}}{4}), (-\frac{5}{2},\frac{\sqrt{3}}{2}), (-\frac{9}{4},\frac{3 \sqrt{3}}{4}), (\frac{5}{2},\frac{\sqrt{3}}{2}), (\frac{11}{4},\frac{\sqrt{3}}{4}), (3,0)\Big\}$ forms an optimal set of seven-points with constrained quantization error $Q_7=\frac{57}{28}$, and there are two different optimal sets of seven-points.
\end{prop1}

 \begin{remark1}
 As the support of $\mu$ contains seven elements, there is no constrained optimal set of $n$-points for any $n\geq 8$.
 \end{remark1}

 \section{Constrained quantization for a finite nonuniform discrete distribution} \label{sec2}
Let $\mu$ be a finite nonuniform distribution with support $\set{-3, -2, -1, 0, 1, 2, 3}$ such that the probability mass function $f$ for $\mu$ is given by
 \[f(j)=\left\{\begin{array}{cc}
 \frac 1{2^{4+j}},  & \te{ where } j\in \set{-3, -2, -1, 0, 1, 2},   \\
 \frac 1{2^6} &\te { if } j=3\\
  0 & \te{otherwise}.
 \end{array}
 \right.\]
In the following two subsections, we calculate the constrained optimal set of $n$-points for the probability distribution $\mu$ for all possible values of $n\in \D N$ with respect to two different constraints.

 \subsection{Constrained quantization when the constraint is a semicircular arc} As mentioned in Subsection~\ref{sub11}, here the constraint is the semicircular arc
 \[R=\set{(3\cos \gq, 3\sin \gq): 0\leq \gq\leq \pi}.\]
Due to the same reasoning as mentioned in Subsection~\ref{sub11}, we can say that the optimal sets of $n$-points exist only for $n=1$ and $n=2$, and they do not exist for any $n\geq 3$.

In the following two propositions, we calculate the optimal sets of $n$-points for $n=1$ and $n=2$.
\begin{prop1}
	The set $\set{(-3, 0)}$ forms the optimal set of one-point with constrained quantization error $Q_1=\frac{177}{64}$.
\end{prop1}

\begin{proof} Let $(3\cos \gq, 3\sin \gq)$ be an element on the circle. Then, the distortion error for $\mu$ with respect to the element is given by
	\begin{align*}
		Q(\mu; \set{(3\cos \gq, 3\sin \gq)})&= \sum_{j=-3}^2 f(j)\rh(j, (3\cos \gq, 3\sin \gq))+f(3)\rh(3,(3\cos \gq, 3\sin \gq))\\
& = \frac{1}{64} (774 \cos \gq +951)
		\end{align*}
	which is minimum when $\gq=\pi$. Hence, the  set $\set{(-3, 0)}$ forms the optimal set of one-point with constrained quantization error $Q_1=\frac{177}{64}$.
\end{proof}

\begin{prop1} \label{prop11}
	The set $\set{(-3, 0), (3, 0)}$ forms an optimal set of two-points with constrained  quantization error $Q_2=\frac{93}{64}$.
\end{prop1}
\begin{proof} Let $\set{(3\cos \gq_1, 3\sin \gq_1), (3\cos \gq_2, 3\sin \gq_2)}$ be an optimal set of two-points, where $0\leq \gq_1<\gq_2\leq \pi$.
As mentioned in Proposition~\ref{prop110}, the support element $0$ falls on the boundary of the two Voronoi regions. Hence, without any loss of generality, we can assume that the support element $0$ belongs to the Voronoi region of $(3\cos \gq_2, 3\sin \gq_2)$. Then, the distortion error is given by
	\begin{align*}
		&Q(\mu; \set{(3\cos \gq_1, 3\sin \gq_1), (3\cos \gq_2, 3\sin \gq_2)})\\
&= \sum_{j=-3}^{0}f(j)\rh(j, (3\cos \gq_2, 3\sin \gq_2))+\sum_{j=1}^{2}f(j)\rh(j, (3\cos \gq_1, 3\sin \gq_1))+f(3)\rh(3, (3\cos \gq_1, 3\sin \gq_1))\\
&=\frac{1}{64} \left(-42 \cos \gq_1+816 \cos \gq_2+951\right)
		\end{align*}
which is minimum if $\gq_1=0$ and $\gq_2=\pi$, and the minimum value is $\frac{93}{64}$. Thus, the set $\set{(-3, 0), (3, 0)}$ forms an optimal set of two-points with constrained  quantization error $Q_2=\frac{93}{64}$.
\end{proof}
\subsection{Constrained quantization when the constraint is the two sides of an equilateral triangle}  As mentioned in Subsection~\ref{sub12}, here the constraint is  $R:=R_1\uu R_2$, where $R_1$ and $R_2$ are the sides $AC$ and $BC$ of the equilateral triangle $ABC$.
 Recall that \begin{align}
R_1=\set{(x, \sqrt{3}(x + 3)) : -3\leq x\leq 0} \te{ and } R_2=\set{(x, \sqrt{3}(-x + 3)) : 0\leq x\leq 3}.
\end{align}
The following proposition gives the main results in this subsection. The techniques of the proof are similar to the proofs given in the propositions in Subsection~\ref{sub12}.
  \begin{prop1} Let $\ga_n$ be an optimal set of $n$-points with constrained quantization error $Q_n$. Then,
   \begin{enumerate}
  \item $\ga_1=\big\{(-\frac{705}{256},\frac{63 \sqrt{3}}{256})\big\}$ with $Q_1=\frac{41343}{16384}$.
      \medskip
    \item $\ga_2=\big\{(-\frac{45}{16},\frac{3 \sqrt{3}}{16}), (\frac{21}{8},\frac{3 \sqrt{3}}{8})\big\}$ with $Q_2=\frac{259}{1024}$.
        \medskip
 \item $\ga_3=\big\{(-\frac{23}{8},\frac{\sqrt{3}}{8}), (-\frac{19}{8},\frac{5 \sqrt{3}}{8}), (\frac{11}{4},\frac{\sqrt{3}}{4})\big\} \te{ with } Q_3=\frac{361}{1536}.$
     \medskip
  \item $\ga_4=\big\{(-3,0), (-\frac{21}{8},\frac{3 \sqrt{3}}{8}), (\frac{19}{8},\frac{5 \sqrt{3}}{8}), (\frac{23}{8},\frac{\sqrt{3}}{8})\big\} \te{ with } Q_4=\frac{117}{512}.$
      \medskip
 \item There are eight different optimal sets of five-points: in five of them three elements are on $R_1$ and two elements are on $R_2$; and in three of them two elements are on $R_1$ and three elements are on $R_2$. One such set is as follows:
     \medskip
 $\ga_5=\big\{(-3,0), (-\frac{21}{8},\frac{3 \sqrt{3}}{8}), (\frac{9}{4},\frac{3 \sqrt{3}}{4}),  (\frac{5}{2},\frac{\sqrt{3}}{2}),  (\frac{23}{8},\frac{\sqrt{3}}{8})\big\} \te{ with } Q_5=\frac{29}{128}$.
     \medskip
 \item There are six different optimal sets of six-points: in five of them three elements are on $R_1$ and three elements are on $R_2$; and in one of them two elements are on $R_1$ and four elements are on $R_2$. One such set is as follows:
     \medskip
 $\ga_6=\big\{(-3,0), (-\frac{21}{8}, \frac{3 \sqrt{3}}{8}), (\frac{9}{4},\frac{3 \sqrt{3}}{4}), (\frac{5}{2},\frac{\sqrt{3}}{2}), (\frac{11}{4},\frac{\sqrt{3}}{4} ),(3,0) \big\} \te{ with } Q_6=\frac{115}{512}$.
     \medskip
 \item $\ga_7=\big\{(-3,0), (-\frac{11}{4},\frac{\sqrt{3}}{4} ), (-\frac{5}{2},\frac{\sqrt{3}}{2}), (-\frac{9}{4},\frac{3 \sqrt{3}}{4}), (\frac{5}{2},\frac{\sqrt{3}}{2} ), (\frac{11}{4},\frac{\sqrt{3}}{4}),(3,0)\big\} \\
     \te{ or } \ga_7=\big\{(-3,0), (-\frac{11}{4},\frac{\sqrt{3}}{4} ), (-\frac{5}{2},\frac{\sqrt{3}}{2}), (\frac{9}{4},\frac{3 \sqrt{3}}{4}), (\frac{5}{2},\frac{\sqrt{3}}{2} ), (\frac{11}{4},\frac{\sqrt{3}}{4}),(3,0)\big\}
      \\ \te{ with } Q_7=\frac{57}{256}.$
  \end{enumerate}
\end{prop1}

%
%
%
%

\section{Constrained quantization for an infinite discrete distribution with support $\set{\frac 1 n : n\in \D N}$} \label{sec3}

Let $\D N:=\set{1, 2, 3, \cdots}$ be the set of natural numbers, and let $\mu$ be a Borel probability measure on the set $\set{\frac 1 n : n \in \D N}$ with probability mass function $f$ given by
\[f(x)=\left\{\begin{array} {ll} \vspace{0.05 in}
\frac 1{2^k} & \te{ if } x =\frac 1 k \te { for } k\in \D N,\\ \vspace{0.05 in}
0 &  \te{ otherwise}.
\end{array}
\right.\]
Then, $\mu$ is a Borel probability measure on $\D R$, and the support of $\mu$ is given by supp$(\mu)=\set{\frac 1 n :  n \in \D N}$.
  In this section, for the probability measure $\mu$, we investigate the constrained optimal sets of $n$-points and the $n$th constrained quantization errors for $n\in \D N$ with respect to two constraints: one is the semicircular arc with base the closed interval $[0, 1]$, and one is the two sides of an equilateral triangle with base $[0, 1]$

 \subsection{Constrained quantization when the constraint is a semicircular arc}
  In this case, let us take the constraint $R$ as the upper semicircular arc of the circle $(x-\frac 12)^2+y^2=\frac 1 4$, i.e.,
 \[R:=\bigg\{(x, y) : \big(x-\frac 12\big)^2+y^2=\frac 12  \te{ and } y\geq 0\bigg\}.\]
 Notice that $R$ can also be represented by
 \[R=\bigg\{ \bigg(\frac 12(1+\cos \gq), \frac 12 \sin \gq\bigg): 0\leq \gq\leq \pi\bigg\}.\]
We know that an optimal set of one-point always exists. For any $n\geq 2$, since the boundary of the Voronoi regions of any two optimal elements, in this case, passes through the center of the circle, from the geometry, we see that among $n$ Voronoi regions, only two Voronoi regions contain elements from the support of $\mu$, i.e., only two Voronoi regions have positive probability. Hence, the optimal sets of $n$-points exist only for $n=1$ and $n=2$, and they do not exist for any $n\geq 3$.

We now calculate the optimal sets of one-point and the two-points in the following propositions.

\begin{prop1}
	An optimal set of one-point is given by $\set{(1, 0)}$  with constrained quantization error $Q_1=\frac{1}{12} \left(\pi ^2-6 \left(-2+\log ^22+\log 16\right)\right)$.
\end{prop1}

\begin{proof} Let $\bigg(\frac 12(1+\cos \gq), \frac 12 \sin \gq\bigg)$ be an element on the circle. Then, the distortion error for $\mu$ with respect to this element is given by
	\begin{align*}
		&Q\Bigg(\mu; \bigg\{\bigg(\frac 12(1+\cos \gq), \frac 12 \sin \gq\bigg)\bigg\} \Bigg)= \sum_{j=1}^\infty \frac 1{2^j} \rh\Bigg(\frac 1 j, \bigg(\frac 12(1+\cos \gq), \frac 12 \sin \gq\bigg)\Bigg) \\
&=\frac{1}{12} \left((6-12 \log 2) \cos \theta +\pi ^2-6 \left(-1+\log ^22+\log 4\right)\right),
		\end{align*}
	which is minimum when $\gq=0$, and the minimum value is $\frac{1}{12} \left(\pi ^2-6 \left(-2+\log ^22+\log 16\right)\right)$. Thus, the proposition is yielded.
\end{proof}

\begin{prop1} \label{prop11}
	The set $\set{(1, 0), (0, 0)}$ forms an optimal set of two-points with constrained  quantization error $Q_2=\frac{1}{12} \left(\pi ^2-6-6 \log ^22-12 \log 2+6 \log 4\right)$.
\end{prop1}
\begin{proof} Let $\ga_2:=\set{(\frac 12 (1+\cos \gq_1), \frac 12 \sin \gq_1), (\frac 12 (1+\cos \gq_2), \frac 12 \sin \gq_2)}$ be an optimal set of two-points, where $0\leq \gq_1<\gq_2\leq \pi$.
Since the boundary of any two elements on the circle passes through the center of the circle, in an optimal set of two-points, one Voronoi region will contain the elements to the left of $\frac 12$, and the other Voronoi region will contain the elements to the right of $\frac 12$. Notice that the support element $\frac 12$ falls on the boundary of the two Voronoi regions. In that case, as mentioned in Remark~\ref{rem1}, we can assume that  $\frac 12$ belongs to the Voronoi region of $(\frac 12 (1+\cos \gq_2), \frac 12 \sin \gq_2)$. Then, the distortion error is given by
	\begin{align*}
		&Q(\mu; \ga)\\
&= \sum_{j=2}^\infty \frac 1{2^j} \rh\Big(\frac 1 j, (\frac 12 (1+\cos \gq_2), \frac 12 \sin \gq_2)\Big)+\frac 12\rh\Big(1, (\frac 12 (1+\cos \gq_1), \frac 12 \sin \gq_1)\Big)\\
&=\frac{1}{12} \left(-3 \cos \theta _1+(9-6 \log 4) \cos \theta _2+\pi ^2+6-6 \log ^22-12 \log 2\right)
		\end{align*}
which is minimum if $\gq_1=0$ and $\gq_2=\pi$, and  $\frac{1}{12} \left(\pi ^2-6-6 \log ^22-12 \log 2+6 \log 4\right)$ is the minimum value. Thus, the proposition is yielded.
\end{proof}
\subsection{Constrained quantization when the constraint is the two sides of an equilateral triangle}

Let us consider the equilateral triangle with vertices $A(0, 0)$, $B(1, 0)$, and $C(\frac 12, \frac{\sqrt 3}2)$. Let $R_1$ and $R_2$ represent the sides $BC$ and $AC$, respectively. In this subsection, we take the constraint as $R:=R_1\uu R_2$. Notice that $R_1$ and $R_2$ can be represented as follows:
\begin{align}
R_1=\big\{(x, \sqrt{3}(1-x)) : \frac 12 \leq x\leq 1\big\} \te{ and }  R_2=\big\{(x, \sqrt{3}x) : 0\leq x\leq \frac 12\big\}.
\end{align}
 For $k, \ell\in \D N$, where $k\leq \ell$, write
\[[k, \ell]:=\Big\{\frac 1 n : n \in \D N \te{ and } k\leq n\leq \ell\Big\}, \te{ and } [k, \infty]:=\Big\{\frac 1 n : n\in \D N \te{ and } n\geq k\Big\}.\]
Further, write
\[Av[k, \ell]: =E\Big (X : X \in [k, \ell]\Big)=\frac{\sum _{n=k}^{\ell} \frac{1}{2^n} \frac 1 n}{\sum_{n=k}^{\ell}\frac{1}{2^n}}, \  Av[k, \infty]: =E\Big (X : X \in [k, \infty)\Big)=\frac{\sum _{n=k}^{\infty} \frac{1}{2^n} \frac 1 n}{\sum_{n=k}^{\infty}\frac{1}{2^n}},\]
\[Er[k, \ell]:=\sum _{n=k}^{\ell} \frac{1}{2^n} \Big(\frac 1 n-Av[k, \ell]\Big)^2, \te{ and } Er[k, \infty]:=\sum _{n=k}^{\infty} \frac{1}{2^n} \Big(\frac 1 n-Av[k,\infty)\Big)^2.\]

\begin{note} \label{note11}
The perpendicular through a point $(a, \sqrt 3(1-a))$ on the line $y=\sqrt 3(1-x)$ intersects the segment $[0, 1]$ at the point $(4a-3,0)$, where  $0\leq 4a-3\leq 1$, i.e., $\frac 34\leq a\leq 1$.
The perpendicular through a point $(b, b\sqrt 3)$ on the line $y= \sqrt 3x$ intersects the segment $[0, 1]$ at the point $(4b, 0)$, where $0\leq 4b\leq 1$, i.e., $0\leq b\leq \frac 14$.
Write
\begin{align*}
R_1^{(1)}:&=\big\{(x, \sqrt 3(1-x)) : \frac 34\leq x\leq 1\big\}, \te{ and } \\
R_2^{(2)}:&=\big\{(x, \sqrt 3x) : 0\leq x\leq \frac 14\big\}.
\end{align*}
Then, there exist two bijective functions $U_1$ and $U_2$, such that

\begin{equation} \label{eq1}
\begin{aligned} &U_1 : [0, 1] \to R_1^{(1)} \te{ such that } U_1(x)=\Big(\frac 14 (3+x), \frac{\sqrt 3}4 (1-x)\Big), \te { and } \\
&U_2 : [0, 1] \to R_2^{(2)} \te{ such that } U_2(x)=\Big(\frac 14 x, \frac{\sqrt 3}4 x\Big).
\end{aligned}
\end{equation}
Let $\ga_n$ be an optimal set of $n$-points for some $n\in \D N$. Then, we must have
\[\ga_n=(\ga_n\ii R_1^{(1)})\uu (\ga_n\ii R_2^{(2)}).\]
Suppose that both $\ga_n\ii R_1^{(1)}$ and $\ga_n\ii R_2^{(2)}$ are nonempty. Then, for some $\ell, m\in \D N$, we can write
\begin{align*}
\ga_n\ii R_1^{(1)} &=\set{(a_\ell,  \sqrt 3(1-a_\ell)), (a_{\ell-1}, \sqrt 3 (1-a_{\ell-1})), \cdots, (a_1, \sqrt 3(1-a_1))}, \te{ and }\\
 \ga_n\ii R_2^{(2)} & =\set{ (a_{\ell+m},  \sqrt 3 a_{\ell+m}), (a_{\ell+m-1},  \sqrt 3 a_{\ell+m-1}), \cdots, (a_{\ell+2},  \sqrt 3 a_{\ell+2}), (a_{\ell+1},  \sqrt 3 a_{\ell+1}))}.
 \end{align*}
 such they satisfy the following order:
 \[ 0\leq a_{\ell+m}<a_{\ell+m-1}<\cdots<a_{\ell+1}\leq \frac 12 \leq a_{\ell} <a_{\ell-1}<\cdots<a_2\leq a_1.\]
Then, the elements $(a_i, \sqrt 3(1-a_i))$ and $(a_{i+1}, \sqrt 3(1-a_{i+1}))$ for $1\leq i\leq \ell-1$;  $(a_\ell,\sqrt 3(1-a_\ell))$ and $(a_{\ell+1},  \sqrt 3 a_{\ell+1})$; and $(a_{\ell+j},  \sqrt 3 a_{\ell+j})$ and $(a_{\ell+j+1},  \sqrt 3 a_{\ell+j+1})$ for $1\leq j\leq m-1$, are the adjacent elements.
\end{note}

\begin{notation}
To be specific, in this subsection,  $Q_n^{(c)}$ represents the constrained quantization error for $n$-points, and $Q_n^{(u)}$ represents the unconstrained quantization error for $n$-means with respect to the probability distribution $\mu$.
\end{notation}

The following proposition plays an important role.

\begin{prop1} \label{propim}
Let $\ga_n:=\set{c_j : 1\leq j\leq k+m}$, where $\te{card}(\ga_n\ii R_1^{(1)})=k$ and $\te{card}(\ga_n\ii R_1^{(2)})=m$. Then,  $\ga_n$ is a  constrained optimal set of $n$-points for $\mu$ iff  $U_1^{-1}(\ga_n\ii R_1^{(1)})\uu U_2^{-1}(\ga_n\ii R_2^{(2)})$ forms an unconstrained optimal set of $n$-means for $\mu$, and \begin{equation*} Q_n^{(c)}=Q_n^{(u)}+\frac 34 \sum_{j=1}^{k}\sum_{i=q_j^{(1)}}^{q_j^{(2)}}\frac 1{2^i}\big(1-Av[q_j^{(1)}, q_j^{(2)}]\big)^2 +\frac 34 \sum_{j=k+1}^{k+m}\sum_{i=q_j^{(1)}}^{q_j^{(2)}}\frac 1{2^i}\big(Av[q_j^{(1)}, q_j^{(2)}]\big)^2,
\end{equation*}
where the Voronoi region of $c_j$ contains the set $[q_j^{(1)}, q_j^{(2)}]$ from the support of $\mu$ for $1\leq j\leq n$, where $q_n^{(2)}=\infty$.
\end{prop1}

\begin{proof}
Let $\ga_n$ be a constrained optimal set of $n$-points for $\mu$. Let $\te{card}(\ga_n\ii R_1^{(1)})=k$ and $\te{card}(\ga_n\ii R_2^{(2)})=m$. Then, we can write
\[\ga_n=\Big\{(a_1, \sqrt 3(1-a_1)), (a_2, \sqrt 3(1-a_2)), \cdots, (a_k, \sqrt 3 (1-a_k)), (a_{k+1}, \sqrt 3a_{k+1}), \cdots, (a_{k+m}, \sqrt 3 a_{k+m})\Big\},\]
where $0\leq a_{k+m}<a_{k+m-1}<\cdots<a_{k+1}\leq \frac 12 \leq a_k<\cdots<a_2<a_1\leq 1$ and $k+m=n$. Let $c_j=(a_j, \sqrt 3(1-a_j))$ if $1\leq j\leq k$, and $c_j=(a_j, \sqrt 3a_j)$ if $k+1\leq j\leq k+m$.   Let the Voronoi region of $c_j$ contains the set $[q_j^{(1)}, q_j^{(2)}]$ from the support of $\mu$ for $1\leq j\leq n$, where $q_n^{(2)}=\infty$.
Then, we have
\begin{align*}
Q_n^{(c)}=\sum_{j=1}^k \sum_{i=q_j^{(1)}}^{q_j^{(2)}}\frac 1{2^i} \rh\Big(\frac 1 i, (a_j, \sqrt 3(1-a_j))\Big)+ \sum_{j=k+1}^{k+m} \sum_{i=q_j^{(1)}}^{q_j^{(2)}}\frac 1{2^i} \rh\Big(\frac 1 i, (a_j, \sqrt 3 a_j)\Big).
\end{align*}
Notice that for $1\leq j\leq k$, the distortion error contributed by $c_j$ on its own Voronoi region is given by
\begin{align*}
 &\sum_{i=q_j^{(1)}}^{q_j^{(2)}}\frac 1{2^i} \rh\Big(\frac 1 i, (a_j, \sqrt 3(1-a_j))\Big)\\
 &=\sum_{i=q_j^{(1)}}^{q_j^{(2)}}\frac 1{2^i}\Big((\frac 1 i -a_j)^2+(\sqrt 3(1-a_j))^2\Big)\\
 &=\sum_{i=q_j^{(1)}}^{q_j^{(2)}}\frac 1{2^i} \Big(\Big((\frac 1 i -\frac 14(3+Av[q_j^{(1)}, q_j^{(2)}]))+(\frac 14(3+Av[q_j^{(1)}, q_j^{(2)}])-a_j)\Big)^2 +(\sqrt 3 (1-a_j))^2\Big)\\
&=\sum_{i=q_j^{(1)}}^{q_j^{(2)}}\frac 1{2^i} \Big( (\frac 1 i -\frac 14(3+Av[q_j^{(1)}, q_j^{(2)}]))^2+2(\frac 1 i -\frac 14(3+Av[q_j^{(1)}, q_j^{(2)}]))(\frac 14(3+Av[q_j^{(1)}, q_j^{(2)}])-a_j)\\
&\qquad \qquad + (\frac 14(3+Av[q_j^{(1)}, q_j^{(2)}])-a_j)^2 +(\sqrt 3 (1-a_j))^2\Big)\\
&=\sum_{i=q_j^{(1)}}^{q_j^{(2)}}\frac 1{2^i} \Big((\frac 1 i -\frac 14(3+Av[q_j^{(1)}, q_j^{(2)}]))^2+(\sqrt 3 (1-a_j))^2\Big) \\
&+2(\frac 14(3+Av[q_j^{(1)}, q_j^{(2)}])-a_j)\sum_{i=q_j^{(1)}}^{q_j^{(2)}}\frac 1{2^i}(\frac 1 i -\frac 14(3+Av[q_j^{(1)}, q_j^{(2)}]))+(\frac 14(3+Av[q_j^{(1)}, q_j^{(2)}])-a_j)^2 \sum_{i=q_j^{(1)}}^{q_j^{(2)}}\frac 1{2^i}  \\
&=\sum_{i=q_j^{(1)}}^{q_j^{(2)}}\frac 1{2^i} \Big((\frac 1 i -\frac 14(3+Av[q_j^{(1)}, q_j^{(2)}]))^2+(\sqrt 3 (1-a_j))^2\Big)  +(\frac 14(3+Av[q_j^{(1)}, q_j^{(2)}])-a_j)^2 \sum_{i=q_j^{(1)}}^{q_j^{(2)}}\frac 1{2^i}.
\end{align*}
Notice that the above distortion error is minimum if $a_j=\frac 14(3+Av[q_j^{(1)}, q_j^{(2)}])$,
 and then
\[c_j=\big(\frac 14(3+Av[q_j^{(1)}, q_j^{(2)}]), \frac{\sqrt 3}4 (1- Av[q_j^{(1)}, q_j^{(2)}])\big)=U_1(Av[q_j^{(1)}, q_j^{(2)}])\] and the distortion error contributed by $c_j$ on its own Voronoi region is given by
\[\sum_{i=q_j^{(1)}}^{q_j^{(2)}}\frac 1{2^i} \rh(\frac 1 i, U_1(Av[q_j^{(1)}, q_j^{(2)}])).\]
Similarly, for $k+1\leq j\leq k+m$, the distortion error  $\sum_{i=q_j^{(1)}}^{q_j^{(2)}}\frac 1{2^i} \rh\Big(\frac 1 i, (a_j, \sqrt 3 a_j)\Big)$ is minimum if $a_j=\frac 14 Av[q_j^{(1)}, q_j^{(2)}]$, and then
\[c_j=\big(\frac 14 Av[q_j^{(1)}, q_j^{(2)}], \frac{\sqrt 3}4 Av[q_j^{(1)}, q_j^{(2)}]\big)=U_2(Av[q_j^{(1)}, q_j^{(2)}])\] and the distortion error contributed by $c_j$ on its own Voronoi region is given by
\[\sum_{i=q_j^{(1)}}^{q_j^{(2)}}\frac 1{2^i} \rh(\frac 1 i, U_2(Av[q_j^{(1)}, q_j^{(2)}])).\]
Thus, we see that
\begin{equation} \label{eq001}
\ga_n=\Big\{U_1(Av[q_j^{(1)}, q_j^{(2)}]) : 1\leq j\leq k\Big\}\UU\Big\{U_2(Av[q_j^{(1)}, q_j^{(2)}]) : k+1\leq j\leq k+m\Big\}
\end{equation}
forms a constrained optimal set of $n$-points for $\mu$ with constrained quantization error
\begin{align*}
Q_n^{(c)}=\sum_{j=1}^k \sum_{i=q_j^{(1)}}^{q_j^{(2)}}\frac 1{2^i} \rh(\frac 1 i, U_1(Av[q_j^{(1)}, q_j^{(2)}]))+\sum_{j=k+1}^{k+m}\sum_{i=q_j^{(1)}}^{q_j^{(2)}}\frac 1{2^i} \rh(\frac 1 i, U_2(Av[q_j^{(1)}, q_j^{(2)}])).
\end{align*}
Let us consider the set $\gb_n$ given by
\begin{equation} \label{eq002}\gb_n=\Big\{Av[q_j^{(1)}, q_j^{(2)}] : 1\leq j\leq k+m\Big\}.\end{equation}
Notice that the unconstrained distortion error for $\mu$ with respect to the set $\gb_n$, denoted by $Q(\mu; \gb_n)$, is defined as
\[Q(\mu; \gb_n)=\sum_{j=1}^{k+m} \sum_{i=q_j^{(1)}}^{q_j^{(2)}}\frac 1{2^i}(\frac 1 i -Av[q_j^{(1)}, q_j^{(2)}])^2.\]
After some calculation, we see that
\begin{equation} \label{eqim} Q_n^{(c)}=Q(\mu; \gb_n)+\frac 34 \sum_{j=1}^{k}\sum_{i=q_j^{(1)}}^{q_j^{(2)}}\frac 1{2^i}(1-Av[q_j^{(1)}, q_j^{(2)}])^2 +\frac 34 \sum_{j=k+1}^{k+m}\sum_{i=q_j^{(1)}}^{q_j^{(2)}}\frac 1{2^i}(Av[q_j^{(1)}, q_j^{(2)}])^2.
\end{equation}
From the above expression, we see that $Q_n^{(c)}$ will be minimum iff the distortion error $Q(\mu; \gb_n)$ is minimum with respect to all possible subsets of $\te{supp}(\mu)$ with cardinality $n$, i.e., when $Q(\mu; \gb_n)$ becomes $Q_n^{(u)}$. Thus, we can say that the set $\ga_n$ given by \eqref{eq001}  will be a constrained optimal set of $n$-points for $\mu$ iff the set $\gb_n$ given by \eqref{eq002}  becomes an unconstrained optimal set of $n$-means for $\mu$.

As $U_1$ and $U_2$ are bijective functions, by \eqref{eq001} and \eqref{eq002}, we see that if $\ga_n$ is a constrained optimal set of $n$-points and $\gb_n$ is an unconstrained optimal set of $n$-means, then
\[U_1^{-1}(\ga_n\ii R_1^{(1)})\uu U_2^{-1}(\ga_n\ii R_2^{(2)})=\gb_n,\]
and by \eqref{eqim}, we have
\begin{equation*} Q_n^{(c)}=Q_n^{(u)}+\frac 34 \sum_{j=1}^{k}\sum_{i=q_j^{(1)}}^{q_j^{(2)}}\frac 1{2^i}(1-Av[q_j^{(1)}, q_j^{(2)}])^2 +\frac 34 \sum_{j=k+1}^{k+m}\sum_{i=q_j^{(1)}}^{q_j^{(2)}}\frac 1{2^i}(Av[q_j^{(1)}, q_j^{(2)}])^2.
\end{equation*} Thus, the proof of the proposition is complete.
\end{proof}
 \begin{prop1}\label{propD1}
An constrained optimal set of one-point is given by $\set{(\frac{1}{4} (3+\log 2),\frac{1}{4} \sqrt{3} (1-\log 2))}$ with constrained quantization error \[Q_1^{(c)}=\frac{1}{12} \left(\pi ^2-9 \left(-1+\log ^22+\log 4\right)\right).\]
\end{prop1}

\begin{proof} Let $\ga:=\set{(a, b)}$ be an optimal set of one-point. The following two cases can happen:

\tit{Case~1. $(a, b)\in R_1$.}

In this case, $b= \sqrt{3}(1-a)$ for some $\frac 12\leq a\leq 1$.  Then, the distortion error is given by
\[Q(\mu; \ga)=\sum_{j=1}^\infty \frac 1{2^j} \rh\big(\frac{1}{j}, (a, b)\big)=4 a^2-2 a (3+\log 2)+\frac{1}{12} \left(\pi ^2+36-6 \log ^22\right),\]
the minimum value of which is $\frac{1}{12} \left(\pi ^2-9 \left(-1+\log ^22+\log 4\right)\right)\approx 0.172407$, and it occurs when $a=\frac{1}{4} (3+\log 2)$.

\tit{Case~2. $(a, b)\in R_2$.}

In this case, $b=\sqrt{3}a$ for some $0\leq a\leq \frac 12$. Then, the distortion error is given by
\[Q(\mu; \ga)=\sum_{j=1}^\infty \frac 1{2^j} \rh\big(\frac{1}{j}, (a, b)\big)=4 a^2-a \log 4+\frac{1}{12} \left(\pi ^2-6 \log ^22\right),\]
the minimum value of which is $\frac{1}{12} \left(\pi ^2-\log ^28\right)\approx 0.462127$, and it occurs when $a=\frac{\log 2}{4}$.

Comparing the two distortion errors, we deduce that an optimal set of one-point is given by $\set{(\frac{1}{4} (3+\log 2),-\frac{1}{4} \sqrt{3} (\log 2-1))}$ with constrained quantization error $Q_1=\frac{1}{12} \left(\pi ^2-9 \left(-1+\log ^22+\log 4\right)\right)$.
\end{proof}

 \begin{remark1} \label{rem11}  Notice that
$E(X):=E(X : X \in \te{supp}(\mu)) =\sum _{n=1}^\infty \frac 1{2^n} \frac 1 n=Av[1, \infty]=\log 2$, and so the unconstrained optimal set of one-mean is the set $\set{\log 2}$ with unconstrained quantization error
\[Q_1^{(u)}=\sum _{n=1}^{\infty} \frac 1 {2^n} \Big (\frac{1}{n}-\log 2\Big)^2=Er[1, \infty] =\frac{1}{12} \left(\pi ^2-18 \log ^2 2\right)=0.101788.\]
By Proposition~\ref{propD1}, we have $\ga_1=\set{(\frac{1}{4} (3+\log 2),\frac{1}{4} \sqrt{3} (1-\log 2))}$. Notice that
\[U_1^{-1}(\ga_1\ii R_1^{(1)})\uu U_2^{-1}(\ga_n\ii R_2^{(2)})=U_1^{-1}(\ga_1\ii R_1^{(1)})=\set{\log 2},\]
which supports Proposition~\ref{propim} for $n=1$.
\end{remark1}
\begin{prop1} \label{propD2}
An optimal set of two-points is given by
$\set{(1,0), (\frac{1}{4} (2\log 2-1),\frac{1}{4} \sqrt{3} (2\log 2-1))}$ with constrained quantization error
\[ Q_2^{(c)}=\frac{\pi ^2}{12}+\frac{1}{8} \left(-5-8 \log ^22+\log 16\right)\approx 0.0635876.\]
\end{prop1}

\begin{proof}

Let $\ga_2$ be an optimal set of two-points.

The following cases can arise:

\tit{Case~1: Both the elements are on $R_1$.}

Then, we can write $\ga_2=\set{(a_1, \sqrt 3(1-a_1)),  (a_2, \sqrt 3(1-a_2))}$ such that $\frac 12\leq a_2<a_1\leq 1$.

Notice that the Voronoi region of $(a_1, \sqrt 3(1-a_1))$ must contain the element $1$ from the support of $\mu$. Assume that the Voronoi region of $(a_1, \sqrt 3(1-a_1))$ contains first $m$ elements from the support of $\mu$. Then, the remaining elements will be contained in the Voronoi region of $(a_2, \sqrt 3(1-a_2))$. Let $F(m)$ denote the minimum distortion error when the Voronoi region of $(a_1, \sqrt 3(1-a_1))$ contains first $m$ elements from the support of $\mu$. Then, we have
\begin{align*}
F(m)=\min\Big\{ \sum_{j=1}^m \frac 1{2^j} \rh\Big(\frac 1 j, (a_1, \sqrt 3(1-a_1))\Big)+\sum_{j=m+1}^\infty \frac 1{2^j} \rh\Big(\frac 1 j, (a_2, \sqrt 3(1-a_2))\Big)   : \frac 12\leq a_2<a_1\leq 1\Big\}.
\end{align*}
Using some simple calculation, we see that $F(m)$ is an increasing sequence in $m\in \D N$. Thus, we see that
\[\min\set{F(m) : m\in \D N}=F(1)=-2+\frac{1}{12} \left(\pi ^2+36-6 \log ^22\right)-\frac{1}{8} (2+\log 4)^2\approx 0.148867,\]
 and it occurs when $a_1=1$ and $a_2=\frac{1}{4} (2+\log 4)$.

\tit{Case~2: One element is on $R_1$ and one element is on $R_2$.}

Then, we can write $\ga_2=\set{(a_1, \sqrt 3(1-a_1)),  (a_2, \sqrt 3 a_2)}$ such that $0 \leq a_2\leq \frac 12 \leq a_1\leq 1$.

Assume that the Voronoi region of $(a_1, \sqrt 3(1-a_1))$ contains first $m$ elements from the support of $\mu$. Then, the remaining elements will be contained in the Voronoi region of $(a_2, \sqrt 3 a_2)$. Let $G(m)$ denote the minimum distortion error when the Voronoi region of $(a_1, \sqrt 3(1-a_1))$ contains first $m$ elements from the support of $\mu$. Then, we have
\begin{align*}
G(m)=\min\Big\{ \sum_{j=1}^m \frac 1{2^j} \rh\Big(\frac 1 j, (a_1, \sqrt 3(1-a_1))\Big)+\sum_{j=m+1}^\infty \frac 1{2^j} \rh\Big(\frac 1 j, (a_2, \sqrt 3 a_2)\Big)   :0 \leq a_2\leq \frac 12 \leq a_1\leq 1\Big\}.
\end{align*}
Using some simple calculation, we see that $G(m)$ is an increasing sequence in $m\in \D N$. Thus, we see that
\[\min\set{G(m) : m\in \D N}=G(1)=\frac{\pi ^2}{12}+\frac{1}{8} \left(-5-8 \log ^22+\log 16\right)\approx 0.0635876,\]
 and it occurs when $a_1=1$ and $a_2=\frac{1}{4} (\log 4-1)$.

\tit{Case~3: Both elements are on $R_2$.}

Then, we can write $\ga_2=\set{(a_1, \sqrt 3 a_1),  (a_2, \sqrt 3 a_2}$ such that $0 \leq a_2< a_1\leq \frac 12$.

Assume that the Voronoi region of $(a_1, \sqrt3 a_1)$ contains first $m$ elements from the support of $\mu$. Then, the remaining elements will be contained in the Voronoi region of $(a_2, \sqrt 3 a_2)$. Let $H(m)$ denote the minimum distortion error when the Voronoi region of $(a_1, \sqrt 3a_1)$ contains first $m$ elements from the support of $\mu$. Then, we have
\begin{align*}
H(m)=\min\Big\{ \sum_{j=1}^m \frac 1{2^j} \rh\Big(\frac 1 j, (a_1, \sqrt 3a_1)\Big)+\sum_{j=m+1}^\infty \frac 1{2^j} \rh\Big(\frac 1 j, (a_2, \sqrt 3 a_2)\Big)   :0 \leq a_2\leq \frac 12 \leq a_1\leq 1\Big\}.
\end{align*}
Using some simple calculation, we see that $H(m)$ is an increasing sequence in $m\in \D N$. Thus, we see that
\[\min\set{H(m) : m\in \D N}=H(1)=\frac{\pi ^2}{12}+\frac{1}{8} \left(-1-8 \log ^2(2)+\log 16\right)\approx 0.563588,\]
 and it occurs when $a_1=\frac 12$ and $a_2=\frac{1}{4} (\log 4-1)$.

Comparing Case~1 through Case~3, we can deduce that an optimal set of two-points occurs in Case~2, i.e.,  an optimal set of two-points is given by
$\set{(1,0), (\frac{1}{4} (\log 4-1),\frac{1}{4} \sqrt{3} (\log 4-1))}$ with constrained quantization error
\[ Q_2^{(c)}=\frac{\pi ^2}{12}+\frac{1}{8} \left(-5-8 \log ^22+\log 16\right)\approx 0.0635876,\]
which yields the proposition.
\end{proof}

\begin{remark1}
By Proposition~\ref{propD2}, we have $\ga_2=\set{(1,0), (\frac{1}{4} (2\log 2-1),\frac{1}{4} \sqrt{3} (2\log 2-1))}$. Then, by Proposition~\ref{propim}, we obtain an unconstrained optimal set of two-means as
\[U_1^{-1}(\ga_2\ii R_1^{(1)})\uu U_2^{-1}(\ga_2\ii R_2^{(2)})=\set{1, 2\log 2-1}.\]
In fact, $\set{1, 2\log 2-1}$ forms an unconstrained optimal set of two-means (see \cite {CHMR}).
\end{remark1}

\begin{lemma1} \label{lemma22}
Let $\ga_n$ be a constrained optimal set of $n$-points for $n\geq 2$. Then, $\ga_n$ contains elements from both $R_1$ and $R_2$.
\end{lemma1}

\begin{proof}
By Proposition~\ref{propD2}, the lemma is true for $n=2$. We now prove that the lemma is true for $n\geq 3$.
Let $\ga_n$ be an optimal set of $n$-points for some $n\geq 3$, and the corresponding $n$th constrained quantization error is $Q_n$. Then,
$Q_n\leq Q_2$, i.e., $Q_n\leq  0.0635876$. For the sake of contradiction, assume that $\ga_n$ does not contain any element from $R_1$. Then, the element $1$ from the support of $\mu$ must be contained in the Voronoi region of an element $(a_1, \sqrt  3 a_1)\in \ga_n\ii R_2$. Then,
\[Q_n\geq \frac 1{2} \rh\Big(1,  (a_1, \sqrt 3 a_1)\Big)=2 a_1^2-a_1+\frac{1}{2},\]
the minimum value of which is $\frac{3}{8}$ and it occurs when $a_1=\frac 14$. Thus, we have $Q_n\geq \frac 38>Q_n$, which is a contradiction.
Hence, $\ga_n$ contains an element from $R_1$. Next, for the sake of contradiction, assume that $\ga_n$ does not contain any element from $R_2$. Then, we can write
\[\ga_n:=\set{(a_1, \sqrt 3(1-a_1)), (a_2, \sqrt 3(1-a_2)), \cdots}, \te{ where } \frac 12\leq \cdots<a_2<a_1.\]
 Clearly the Voronoi region of $(a_1, \sqrt 3 (1-a_1))$ contains $1$ from the support of $\mu$. Suppose that the Voronoi region of $(a_1, \sqrt 3 (1-a_1))$ also contains $1$, $\frac 12$ and $\frac 13$. Then,
\[Q_n\geq \sum_{j=1}^3 \frac 1{2^j} \rh\Big(\frac 1 j,  (a_1, \sqrt 3(1-a_1))\Big)=\frac{1}{144} \left(12 a_1 \left(42 a_1-79\right)+461\right),\]
the minimum value of which is $\frac{71}{672}$ and it occurs when $a_1=\frac{79}{84}$. Thus, we see that $Q_n\geq \frac{79}{84}>Q_n$, which is a contradiction.
 Suppose that the Voronoi region of $(a_1, \sqrt 3 (1-a_1))$ contains only $1$ and $\frac 12$. Then, the Voronoi region of $(a_2, \sqrt 3 (1-a_2))$ must contain $\frac 13$. Then,
\begin{align*}
Q_n&\geq \sum_{j=1}^2 \frac 1{2^j} \rh\Big(\frac 1 j,  (a_1, \sqrt 3(1-a_1))\Big)+\frac 1{2^3} \rh\Big(\frac 13,  (a_2, \sqrt 3(1-a_2))\Big)=3 a_1^2-\frac{23 a_1}{4}+\frac{1}{6} a_2 \left(3 a_2-5\right)+\frac{461}{144},
\end{align*}
the minimum value of which is $\frac{19}{192}$ and it occurs when $a_1=\frac{23}{24}$ and $a_2=\frac{5}{6}$. Thus, we see that $Q_n\geq \frac{19}{192}>Q_n$, which leads to a contradiction. Hence, we can assume that the Voronoi region of $(a_1, \sqrt 3 (1-a_1))$ contains only the element $1$. Hence, the Voronoi region of $(a_2, \sqrt 3 (1-a_2))$ must contain the element $\frac 12$ from the support of $\mu$. Suppose that it also contains $\frac 13$. Then,
\begin{align*}
Q_n&\geq   \frac 1{2} \rh\Big(1,  (a_1, \sqrt 3(1-a_1))\Big)+\sum_{j=2}^3 \frac 1{2^j} \rh\Big(\frac 1 j,  (a_2, \sqrt 3(1-a_2))\Big)=2 a_1^2-4 a_1+\frac{3 a_2^2}{2}-\frac{31 a_2}{12}+\frac{461}{144},
\end{align*}
the minimum value of which is $\frac{77}{864}$ and it occurs when $a_1=1$ and $a_2= \frac{31}{36}$. Thus, we see that $Q_n\geq \frac{77}{864}>Q_n$, which yields a contradiction.
Hence, we can assume that the Voronoi region of $(a_2, \sqrt 3 (1-a_2))$ contains only the element $\frac 12$. Hence, the Voronoi region of $(a_3, \sqrt 3 (1-a_3))$ must contain the element $\frac 13$. Then,

\begin{align*}
Q_n&\geq   \frac 1{2} \rh\Big(1,  (a_1, \sqrt 3(1-a_1))\Big)+\frac 1{2^2} \rh\Big(\frac 12,  (a_2, \sqrt 3(1-a_2))\Big)+ \frac 1{2^3} \rh\Big(\frac 1 j,  (a_3, \sqrt 3(1-a_3))\Big)\\
&=2 a_1^2-4 a_1+a_2^2+\frac{a_3^2}{2}-\frac{7 a_2}{4}-\frac{5 a_3}{6}+\frac{461}{144},
\end{align*}
the minimum value of which is $\frac{17}{192}$, and it occurs when $a_1= 1, \, a_2=\frac{7}{8}$, and $a_3=\frac{5}{6}$.
Hence, $Q_n\geq \frac{17}{192}>Q_n$, which is a contradiction.

Thus, we see that if $\ga_n$ does not contain any element from $R_2$, we get into a contradiction. Thus, we see that $\ga_n$ contains an element from $R_2$.
In the beginning we have shown that $\ga_n$ contains an element from $R_1$. Thus, the proof of the lemma is complete.
\end{proof}

Recall Lemma~\ref{lemma22}. Then, using a similar technique as Proposition~\ref{propD2}, the following proposition can be proved.

\begin{prop1} \label{propD3}
Let $\ga_3$ be an optimal set of three-points with constrained quantization error $Q_3$. Then,
\begin{align*}
\ga_3=\Big\{(1,0), & (\frac{7}{8},\frac{\sqrt{3}}{8}),  (\frac{1}{8} (8 \log 2-5),\frac{1}{8} \sqrt{3} (8 \log 2-5))\Big\},  \te{ or } \\
& \Big\{(1,0),  (\frac{1}{8},\frac{\sqrt{3}}{8}),  (\frac{1}{8} (8 \log 2-5),\frac{1}{8} \sqrt{3} (8 \log 2-5))\Big\}
\end{align*} with
\[ Q_3^{(c)}=\frac{1}{96} \left(8 \pi ^2-87-144 \log ^22+120 \log 2\right)\approx 0.0619715.\]
\end{prop1}
\begin{remark1}
By Proposition~\ref{propD3}, we know the constrained optimal set $\ga_3$ of three-points. Then, by Proposition~\ref{propim}, we obtain an unconstrained optimal set of three-means as
\[U_1^{-1}(\ga_3\ii R_1^{(1)})\uu U_2^{-1}(\ga_3\ii R_2^{(2)})=\set{1, \frac 12, \frac{1}{2} (8 \log 2-5)}.\]
In fact, $\set{1, \frac 12, \frac{1}{2} (8 \log 2-5)}$ forms an unconstrained optimal set of three-means (see \cite {CHMR}).
\end{remark1}
\begin{lemma1} \label{lemmaIm}
For $ n \geq 3 $, let $ \ga_n $ denote a constrained optimal set of $ n $-points for the probability measure $ \mu $. Then, the cardinality of $ \ga_n \cap R_1 $, denoted by $ \mathrm{card}(\ga_n \cap R_1) $, is either one or two, and not more than two.
\end{lemma1}
\begin{proof}
Let $\ga_n$ be a constrained optimal set of $n$-points for any $n\geq 3$ with constrained quantization error $Q_n$. Then,
\[Q_n\leq Q_3\leq 0.0619715.\]
 By Lemma~\ref{lemma22}, we see that $\te{card}(\ga_n\ii R_1)\geq 1$. We now show that $\te{card}(\ga_n\ii R_1)<3$. For the sake of contradiction assume that $\te{card}(\ga_n\ii R_1)\geq 3$. Then,  we can write
\[\ga_n:=\set{(a_1, \sqrt 3(1-a_1)), (a_2, \sqrt 3(1-a_2)), \cdots}, \te{ where } \frac 12\leq \cdots<a_2<a_1.\]
As shown in Lemma~\ref{lemma22}, we can see that the Voronoi region of $(a_1, \sqrt 3 (1-a_1))$ contains only the element $1$, and the Voronoi region of $(a_2, \sqrt 3 (1-a_2))$ contains only the element $\frac 12$. Hence, the Voronoi region of $(a_3, \sqrt 3 (1-a_3))$ must contain the element $\frac 13$. Then,
\begin{align*}
Q_n&\geq   \frac 1{2} \rh\Big(1,  (a_1, \sqrt 3(1-a_1))\Big)+\frac 1{2^2} \rh\Big(\frac 12,  (a_2, \sqrt 3(1-a_2))\Big)+ \frac 1{2^3} \rh\Big(\frac 1 j,  (a_3, \sqrt 3(1-a_3))\Big)\\
&=2 a_1^2-4 a_1+a_2^2+\frac{a_3^2}{2}-\frac{7 a_2}{4}-\frac{5 a_3}{6}+\frac{461}{144},
\end{align*}
the minimum value of which is $\frac{17}{192}$, and it occurs when $a_1= 1, \, a_2=\frac{7}{8}$, and $a_3=\frac{5}{6}$.
Hence, $Q_n\geq \frac{17}{192}>Q_n$, which is a contradiction. Hence, $\te{card}(\ga_n\ii R_1)<3$. We now show that $\te{card}(\ga_n\ii R_1)$ can be either one or two. By Proposition~\ref{propD3}, we see that $\te{card}(\ga_n\ii R_1)$ can be either one or two for $n=3$. We now claim that $\te{card}(\ga_n\ii R_1)$ can be either one or two for $n\geq 4$. The distortion error contributed by the set
\[\gb:= \Big\{(1,0), (\frac{1}{8},\frac{\sqrt{3}}{8}), (\frac{1}{12},\frac{1}{4 \sqrt{3}}),  (\frac{2}{3} (3 \log (2)-2),\frac{2 (3 \log (2)-2)}{\sqrt{3}})\Big\},\]
is given by
\begin{align*}
 Q(\mu; \gb)&=\frac 1{2} \rh\Big(1,  (1,0))\Big)+ \frac 1{2^2} \rh\Big(\frac 1 2,  (\frac{1}{8},\frac{\sqrt{3}}{8})\Big)+ \frac 1{2^3} \rh\Big(\frac 1 2,  (\frac{1}{12},\frac{1}{4 \sqrt{3}})\Big)\\
 & +\sum_{j=4}^\infty \frac 1{2^j} \rh\Big(\frac 1 j,   (\frac{2}{3} (3 \log 2-2),\frac{2 (3 \log 2-2)}{\sqrt{3}})\Big)=0.0617409.
\end{align*}
Let $\ga_n:=\set{c_1, c_2, c_3, \cdots}$ be a constrained optimal set of $n$-points for $n\geq 4$ with constrained quantization error $Q_n$. Then,
\[Q_n\leq Q_4\leq0.0617409.\] By Lemma~\ref{lemma22}, we know $c_1\in R_1$.
For sure, we know that $c_j\in R_2$ for all $j\geq 3$. We now show that the Voronoi region of $c_2$ contains only the element $\frac 12$.

The following two cases can arise.

\tit{Case~1. $c_2\in R_1$.}

Then, as shown in the proof of Lemma~\ref{lemma22}, we can say that the Voronoi region of $c_2$ contains only the element $\frac 12$.

\tit{Case~2. $c_2\in R_2$.}

Then, we can write $c_1=(a_1, \sqrt 3(1-a_1))$ and $c_j=(a_j, \sqrt 3 a_j)$ for $j\geq 2$ such that $0<\cdots <a_{j+1}<a_j\leq \frac 12 \leq  a_1\leq 1$.  The Voronoi region of $(a_1, \sqrt 3(1-a_1))$ must contain $1$. Suppose that it also contains $\frac 12$ and $\frac 13$.
Then,
\[Q_n\geq \sum_{j=1}^3\frac 1 {2^j} \rh\Big(\frac 1 j,  (a_1, \sqrt 3(1-a_1))\Big)\geq \frac{71}{672}>Q_n,\]
and the minimum value occurs when $a_1=\frac{79}{84}$.
 Thus, we arrive at a contradiction. Now, assume that the Voronoi region of $(a_1, \sqrt 3(1-a_1))$ contains only the elements $1$ and $\frac 12$. Then, the Voronoi region of $(a_2, \sqrt 3 a_2)$ must contain $\frac 13$. Then,
 \[Q_n\geq \sum_{j=1}^2\frac 1 {2^j} \rh\Big(\frac 1 j,  (a_1, \sqrt 3(1-a_1))\Big)+\frac 1{2^3}  \rh\Big(\frac 1 3,  (a_2, \sqrt 3a_2)\Big) \geq\frac{13}{192}>Q_n,\]
 and the minimum value occurs when $a_1=\frac{23}{24}$ and $a_2=\frac 1 {12}$. Thus, we see a contradiction arises. Hence, we can assume that the Voronoi region of $(a_1, \sqrt 3(1-a_1))$ contains only the element $1$. Then, the Voronoi region of $(a_2,\sqrt 3 a_2)$ must contain the element $\frac 12$. Suppose that the Voronoi region of $(a_2,\sqrt 3 a_2)$ also contains the elements $\frac 13$, $\frac 14$ and $\frac 15$. Then,
 \[Q_n\geq \sum_{j=1}^1\frac 1 {2^j} \rh\Big(\frac 1 j,  (a_1, \sqrt 3(1-a_1))\Big)+\sum_{j=2}^5\frac 1 {2^j} \rh\Big(\frac 1 j,   (a_2, \sqrt 3a_2)\Big) \geq\frac{108149}{1728000}>Q_n,\]
 and the minimum value occurs when $a_1=1$ and $a_2=\frac{181}{1800}$. Thus, we get a contradiction arises. Hence, we can assume that the Voronoi region of $(a_2, \sqrt 3 a_2)$ contains only the element $\frac 12$, $\frac 13$ and $\frac 14$. Then,  the Voronoi region of $(a_3, \sqrt 3 a_3)$ must contain the element $\frac 15$. Then,
 \[Q_n\geq \sum_{j=1}^1\frac 1 {2^j} \rh\Big(\frac 1 j,  (a_1, \sqrt 3(1-a_1))\Big)+\sum_{j=2}^4\frac 1 {2^j} \rh\Big(\frac 1 j,   (a_2, \sqrt 3a_2)\Big) +\sum_{j=5}^5\frac 1 {2^j} \rh\Big(\frac 1 j,   (a_3, \sqrt 3a_3)\Big)\geq\frac{14341}{230400}>Q_n,\]
 and the minimum value occurs when $a_1=1$, $a_2=\frac{5}{48}$, and $a_3=\frac{1}{20}$, which is a contradiction. Hence, we can assume that the Voronoi region of $(a_2, \sqrt 3 a_2)$ contains only the element $\frac 12$ and $\frac 13$. Then,  the Voronoi region of $(a_3, \sqrt 3 a_3)$ must contain the element $\frac 14$. Suppose it also contains $\frac 15$.  Then,
 \[Q_n\geq \sum_{j=1}^1\frac 1 {2^j} \rh\Big(\frac 1 j,  (a_1, \sqrt 3(1-a_1))\Big)+\sum_{j=2}^3\frac 1 {2^j} \rh\Big(\frac 1 j,   (a_2, \sqrt 3a_2)\Big) +\sum_{j=4}^5\frac 1 {2^j} \rh\Big(\frac 1 j,   (a_3, \sqrt 3a_3)\Big)\geq \frac{21341}{345600}>Q_n,\]
 and the minimum value occurs when $a_1=1$, $a_2=\frac{1}{9}$, and $a_3=\frac{7}{120}$,  which is a contradiction. Thus, we can assume that  the Voronoi region of $(a_3, \sqrt 3 a_3)$ contains only the element $\frac 14$. Then, the Voronoi region of $(a_4, \sqrt 3 a_4)$ must contain $\frac 15$. Notice that if $n=4$, then the Voronoi region of $(a_4, \sqrt 3 a_4)$ will also contain $\frac 16$, and as before we can see that a contradiction arises. Hence, for $n=4$, we can assume that the Voronoi region of $a_2$ contains only the element $\frac 12$. Suppose that $n\geq 5$.
 Then, if the Voronoi region of $(a_4, \sqrt 3 a_4)$ also contains $\frac 16$, likewise, we can show a contradiction arises. Suppose that the Voronoi region of $(a_4, \sqrt 3 a_4)$ contains only the element $\frac 15$. Then, the Voronoi region of $(a_5, \sqrt 3 a_5)$ must contains $\frac 16$. Then,
\begin{align*}
Q_n&\geq \sum_{j=1}^1\frac 1 {2^j} \rh\Big(\frac 1 j,  (a_1, \sqrt 3(1-a_1))\Big)+\sum_{j=2}^3\frac 1 {2^j} \rh\Big(\frac 1 j,   (a_2, \sqrt 3a_2)\Big) +\sum_{j=4}^4\frac 1 {2^j} \rh\Big(\frac 1 j,   (a_3, \sqrt 3a_3)\Big)\\
& +\sum_{j=5}^5\frac 1 {2^j} \rh\Big(\frac 1 j,   (a_4, \sqrt 3a_4)\Big)+\sum_{j=6}^6\frac 1 {2^j} \rh\Big(\frac 1 j,   (a_5, \sqrt 3a_5)\Big)\geq \frac{21449}{345600}>Q_n,
\end{align*}
 and the minimum value occurs when $a_1=1$, $a_2=\frac 1 9$, $a_3=\frac{1}{16}$, $a_4=\frac{1}{20}$, and $a_5=\frac{1}{24}$, which gives a contradiction. Hence, we can assume that the Voronoi region of $c_2:=(a_2, \sqrt 3 a_2)$ contains only the element $\frac 12$.

 Thus, in both the cases we see that the Voronoi region of $c_2$ contains only the element $\frac 12$. Notice that if $c_2:=(a_2, \sqrt 3(1-a_2))\in R_1$, then the distortion error contributed by $c_2$ in its own Voronoi region is given by
 \[\frac{1}{2^2} \rh\Big(\frac 12,(a_2, \sqrt 3(1-a_2)\Big),\]
the minimum value of which is $\frac 3{64}$, and it occurs when $a_2=\frac 78$. If $c_2:=(a_2, \sqrt 3a_2)\in R_2$, then the distortion error contributed by $c_2$ in its own Voronoi region is given by
 \[\frac{1}{2^2} \rh\Big(\frac 12,(a_2, \sqrt 3a_2\Big),\]
the minimum value of which is $\frac 3{64}$, and it occurs when $a_2=\frac 18$.
Since the distortion error given by $c_2$ when $c_2\in R_1$ is same as the distortion error when $c_2\in R_2$, we can say that in an optimal set of $n$-points for $n\geq 3$, $c_2$ can belong to either $R_1$ or $R_2$, which yields the fact that $\te{card}(\ga_n \cap R_1)$, attain either one or two. We have also shown that $\te{card}(\ga_n \cap R_1)<3$. Thus, the proof of the lemma is complete.
\end{proof}

 \begin{remark1} \label{rem45}
 Because of Lemma~\ref{lemmaIm}, in the sequel in calculating the constrained optimal sets of $n$-points for any $n\geq 4$, we will assume that $\ga_n$ takes only one element from $R_1$ and the remaining $(n-1)$ elements it takes from $R_2$. Moreover, not directly calculating the constrained optimal sets of $n$-points for any $n\geq 4$, we will use our knowledge of unconstrained optimal sets of $n$-means for $n\geq 4$, and then will use Proposition~\ref{propim} to calculate the constrained optimal sets of $n$-points and the corresponding $n$th constrained quantization error as follows:
Let  $\gb_n:=\Big\{Av[q_j^{(1)}, q_j^{(2)}] : 1\leq j\leq n\Big\}$ be an unconstrained optimal set of $n$-means for $\mu$ with the $n$th unconstrained quantization error $Q_n^{(u)}(\mu)$, and let $\ga_n$ be a constrained optimal set of $n$-points for $\mu$ with $n$th constrained quantization error $Q_n^{(c)}(\mu)$. Then,
\begin{equation} \label{eqM1} \begin{aligned}
\ga_n &=\set{U_1(Av[q_1^{(1)}, q_1^{(2)}])}\UU \set{U_2(Av[q_j^{(1)}, q_j^{(2)}]) : 2\leq j\leq n}, \te{ and } \\
Q_n^{(c)}(\mu)&=Q_n^{(u)}(\mu)+ \frac 34 \sum_{j=2}^n\sum_{i=q_j^{(1)}}^{q_j^{(2)}}\frac 1{2^i}(Av[q_j^{(1)}, q_j^{(2)}])^2.
\end{aligned}
\end{equation}
 \end{remark1}
The following theorem gives the constrained optimal set of $n$-points and the corresponding $n$th constrained quantization errors for $1\leq n\leq 2000$.
\begin{theo1} \label{theo1}
For any positive integer $n$, if $1\leq n\leq 5$, then the sets
\[\Big\{U_2(Av[n, \infty)), U_2(\frac 1 {n-1}),  \cdots, U_2(\frac 13), U_2(\frac 12), U_1(1)\Big\},\] where $1\leq n\leq 5$, form  constrained optimal sets of $n$-points;  on the other hand, if $6\leq n\leq 2000$, then the sets $\set{U_2(Av[n+1, \infty)), U_2(Av[n-1, n]), U_2(\frac 1 {n-2}), \cdots, U_2(\frac 13), U_2(\frac 12), U_1(1)}$ form constrained optimal sets of $n$-points for the probability measure $\mu$. The corresponding $n$th constrained quantization errors are given by
\begin{align*}
 Q_n^{(c)}(\mu)=\left\{\begin{array}{cc}
 Er[n, \infty]+\frac 34 \mathop{\sum}\limits_{j=2}^{n-1}\frac 1 {2^j} (\frac 1 j)^2+\frac 34 \mathop{\sum}\limits_{j=n}^\infty \frac 1{2^j} (Av[n, \infty])^2 & \te{ if } 1\leq n\leq 5, \\
\Big(Er[n+1, \infty]+Er[n-1, n]+\frac 34 \mathop{\sum}\limits_{j=2}^{n-2}\frac 1 {2^j} (\frac 1 j)^2+\frac 34\frac 1{2^{n-1}} (Av[n-1, n])^2 \\
+\frac 34 \mathop{\sum}\limits_{j=n}^\infty \frac 1{2^j} (Av[n, \infty])^2\Big)  & \te{ if } 6\leq n\leq 2000.
\end{array}
\right.
\end{align*}
\end{theo1}

 \begin{proof}
 By \cite[Theorem~4.8]{CHMR}, it is known that
for any positive integer $n$, the sets
\begin{equation} \label{eqM4567} \Big\{Av[n, \infty], \frac 1 {n-1},  \cdots, \frac 13, \frac 12, 1\Big\},
\end{equation}  where $1\leq n\leq 5$, form unconstrained optimal sets of $n$-means for the probability measure $\mu$ with unconstrained quantization errors
\begin{equation} \label{eqM4568 }Q_n(\mu):=Er[n, \infty].
\end{equation}
Write
\begin{equation} \label{eqM45}  \ga_n:=\set{Av[n+1, \infty], Av[n-1, n], \frac 1 {n-2}, \cdots, \frac 13, \frac 12, 1}
\end{equation}  and
\begin{equation} \label{eqM46} Q_n(\mu):=Er[n+1, \infty]+Er[n-1, n].
\end{equation}
We now prove the following claim:

\tit{Claim. The set $\ga_n$ given by \eqref{eqM45} for $n=2000$ form an optimal set of $n$-means with unconstrained quantization error $Q_n(\mu)$ given by \eqref{eqM46}. }

The distortion error contributed by the set $\gb:=\set{Av[2001, \infty], Av[1999, 2000], \frac 1{1998}, \cdots, \frac 12, 1}$ is given by
\[Q(\mu; \gb)=Er[2001, \infty]+Er[1999, 2000]=1.4444190239372496100\times 10^{-615}. \]
Since $Q_{2000}$ is the unconstrained quantization error, we have
\[Q_{2000}\leq 1.4444190239372496100\times 10^{-615}.\]
Let $\ga_n:=\set{a_{2000}, a_{1999}, a_{1997}, \cdots, a_2, a_1}$ form an optimal set of $n$-means for $\mu$ with unconstrained quantization error $Q_n$. We show that $a_1=1$. Clearly, the Voronoi region of $a_1$ contains $1$. Suppose it also contains the element $\frac 12$. Then,
\[Q_n\geq Er[1,2]>Q_n,\]
which is a contradiction. Hence, the Voronoi region of $a_1$ contains only the element $1$, i.e., $a_1=1$. Similarly, we can show it $a_i=\frac 1 i$ for all $1\leq i\leq 1997$. Now, we show that $a_{1998}=\frac 1{1998}$. As $a_i=\frac 1i$ for $1\leq i\leq 1997$, the Voronoi region of $a_{1998}$ must contain $\frac 1{1998}$. Suppose that the Voronoi region of $a_{1998}$ also contains $\frac 1{1999}$ and $\frac 1{2000}$. Then,
\[Q_n\geq Er[1998, 2000]=2.0264380104700818609 \times 10^{-615}>Q_n,\]
which leads to a contradiction. Suppose that the Voronoi region of $a_{1998}$ contains only $\frac 1{1998}$ and $\frac 1{1999}$. Then, the Voronoi region of $a_{1999}$ must contain $\frac 1{2000}$. Suppose that the Voronoi region of $a_{1999}$ also contains $\frac 1{i}$ for $i=2001, 2002, 2003$. Then,
\[Q_n\geq Er[1998, 1999]+Er[2000, 2003]=1.6059068732253525792\times 10^{-615}>Q_n,\]
which leads to a contradiction. Suppose that the Voronoi region of $a_{1999}$ contains only the elements $\frac 1{i}$ for $i=2000, 2001, 2002$. Then, the Voronoi region of $a_{2000}$ contains all the remaining elements $\frac 1 i$ for $i\geq 2003$. Then,
\[Q_n\geq Er[1998, 1999]+Er[2000, 2002]+Er[2003, 2010]=1.4671459164463424875\times 10^{-615}>Q_n,\]
which yields a contradiction. Suppose that the Voronoi region of $a_{1999}$ contains only the elements $\frac 1{i}$ for $i=2000, 2001$. Then, the Voronoi region of $a_{2000}$ contains all the remaining elements $\frac 1 i$ for $i\geq 2002$. Then,
\[Q_n\geq Er[1998, 1999]+Er[2000, 2001]+Er[2002, 2015]=1.4455437344001957110\times 10^{-615}>Q_n,\]
which yields a contradiction. Suppose that the Voronoi region of $a_{1999}$ contains only the elements $\frac 1{i}$ for $i=2000$. Then, the Voronoi region of $a_{2000}$ contains all the remaining elements $\frac 1 i$ for $i\geq 2001$. Thus,
\[Q_n\geq Er[1998, 1999]+Er[2001, 2015]=1.8054463072890087690\times 10^{-615}>Q_n,\]
which leads to a contradiction.

Thus, we can assume that the Voronoi region of $a_{1998}$ contains only $\frac 1{1998}$, i.e., $a_{1998}=\frac 1{1998}$. We now show that the Voronoi region of $a_{1999}$ contains only the element $Av[1999, 2000]$. The Voronoi region of $a_{1999}$ must contain $\frac 1{1999}$. Suppose that it also contains $\frac 1 i$ for $i=2000, 2001, 2002$. Then,
\[Q_n\geq Er[1999, 2002]=1.7593293901236548092 \times 10^{-615}>Q_n,\]
which gives  a contradiction. Suppose that the Voronoi region of $a_{1999}$ contains only $\frac 1{1999}$, $\frac 1{2000}$, and $\frac 1{2001}$. Then, the Voronoi region of $a_{2000}$ must contain the remaining elements. Then,
\[Q_n\geq Er[1999, 2001]+Er[2002, 2010]=1.5071509625668501342\times 10^{-615}>Q_n,\]
which yields a contradiction. Suppose that the Voronoi region of $a_{1999}$ contains only $\frac 1{1999}$. Then, the Voronoi region of $a_{2000}$ must contain the remaining elements. Then,
\[Q_n\geq Er[2000, 2010]=1.5159843420767176934\times 10^{-615}>Q_n,\]
which leads to a contradiction. Hene, we can assume that the Voronoi region of $a_{1999}$ contains only the elements $\frac 1{1999}$ and $\frac 1{2000}$, and then the remaining elements are contained in the Voronoi region of $a_{2000}$. Thus, we have
\[a_{1999}=Av[1999, 2000] \te{ and } a_{2000}=Av[2001, \infty],\]
i.e., the set $\ga_n$ given by \eqref{eqM45} for $n=2000$ form an optimal set of $n$-means with unconstrained quantization error $Q_n(\mu)$ given by \eqref{eqM46}. Thus, the claim is true.

Proceeding in the similar way as the proof of the claim, it can be proved that the set $\ga_n$ given by \eqref{eqM45} form optimal sets of $n$-means with unconstrained quantization errors $Q_n(\mu)$ given by \eqref{eqM46} for all $6\leq n\leq 1999$, i.e., the set $\ga_n$ given by \eqref{eqM45} form optimal sets of $n$-means with unconstrained quantization errors $Q_n(\mu)$ given by \eqref{eqM46} for all $6\leq n\leq 2000$. Hence, by Proposition~\ref{propim}, the proof of the theorem is obtained.
 \end{proof}

 As a consequence of Theorem~\ref{theo1}, the following open problem arises.
 \begin{open}
It is still not known whether the sets $\set{Av[n+1, \infty], Av[n-1, n], \frac 1 {n-2}, \cdots, \frac 13, \frac 12, 1}$ give the unconstrained optimal sets of $n$-means for all positive integers $n\geq 2001$. If not, then the least upper bound of $n\in \D N$ for which such sets give the unconstrained optimal sets of $n$-means for the probability measure $\mu$ is not known yet. Consequently, constrained optimal sets of $n$-points for all positive integers $n\geq 2001$ are not known yet.
 \end{open}

\section{Constrained quantization for an infinite discrete distribution with support $\set{n : n\in \D N}$} \label{sec4}

 Let $L$ be a line given by $y=mx+c$ for $c,m\in\D R$, the parametric representation of which is
$L:=\set{(x, mx+c) : x\in \D R}.$
Write
\[
\mu: = \sum_{n=1}^\infty \frac 1 {2^n} \, \delta_n,
\]
where \( \delta_n \) is the \tit{Dirac measure} (also called the \tit{Dirac delta measure}) centered at the point \( n \). That is, for any measurable set \( A \subset \mathbb{R} \),
\[
\delta_n(A) =
\begin{cases}
1 & \text{if } n \in A, \\
0 & \text{if } n \notin A.
\end{cases}
\]
In other words, \( \delta_n \) is a \emph{unit point mass} concentrated at \( n \). Then, $\mu$ is a Borel probability measure on $\D R$ with probability mass function $f$ given by

\[f(x)=\left\{\begin{array}{ccc}
	\frac 1 {2^n} & \te{ if } x=n \te{ for } n \in\D N,\\
	0  & \te{ otherwise}.
\end{array}\right.
\]
In this section, we determine the constrained optimal sets of $n$-points and the $n$th constrained quantization errors for the probability measure $\mu$ for all positive integers $n$ so that the elements in the optimal sets lie on the line $L$, i.e., $L$ plays as a constraint in the constrained quantization.
 For $k, \ell\in \D N$, where $k\leq \ell$, write
\[[k, \ell]:=\set{n : n \in \D N \te{ and } k\leq n\leq \ell}, \te{ and } [k, \infty]:=\set{n : n\in \D N \te{ and } n\geq k}.\]
Further, write
\[Av[k, \ell]: =E\Big (X : X \in [k, \ell]\Big)=\frac{\sum _{n=k}^{\ell} \frac{1}{2^n} n}{\sum_{n=k}^{\ell}\frac{1}{2^n}}, \  Av[k, \infty]: =E\Big (X : X \in [k, \infty)\Big)=\frac{\sum _{n=k}^{\infty} \frac{1}{2^n} n}{\sum_{n=k}^{\infty}\frac{1}{2^n}},\]
\[Er[k, \ell]:=\sum _{n=k}^{\ell} \frac{1}{2^n} \Big(n-Av[k, \ell]\Big)^2, \te{ and } Er[k, \infty]:=\sum _{n=k}^{\infty} \frac{1}{2^n} \Big(n-Av[k,\infty)\Big)^2.\]

	\begin{note1} \label{note11}
The perpendicular through a point $(a, ma+c)$ on the line $y=mx+c$ intersects the real line at the point $(a(1+m^2)+cm, 0)$.
 Define a function
 \[U : \D R \to L \te{ such that } U(x)=\Big(\frac{x-cm}{1+m^2}, \frac{mx+c}{1+m^2}\Big).\]
 Then, $U$ is a bijective function, and its inverse is given by
 \[U^{-1} : L \to \D R \te{ such that } U^{-1}(x, y)=x(1+m^2)+cm.\]
\end{note1}

\begin{prop} \label{propim1}
Let $\ga_n$ be an unconstrained optimal set of $n$-means for $\mu$. Then, $U(\ga_n)$ is a constrained optimal set of $n$-points for $\mu$ with respect to the constraint $L$. The converse is also true, i.e., if $\ga_n$ is a constrained optimal set of $n$-points for $\mu$ with respect to the constraint $L$, then $U^{-1}(\ga_n)$ forms an unconstrained optimal set of $n$-means for $\mu$. Further,
\begin{equation} \label{Megha42} Q_n^{(c)}=Q_n^{(u)}+\frac 1 {1+m^2}\sum_{j=1}^{n}(m Av[q_j^{(1)}, q_j^{(2)}] +c)^2(2^{1-q_j^{(1)}}-2^{-q_j^{(2)}}),
\end{equation}
if $\ga_n:=\set{Av[q_j^{(1)}, q_j^{(2)}] : 1\leq j\leq n}$, where $q_1^{(1)}=1$, $q_n^{(2)}=\infty$ and $q_j^{(2)}+1=q_{j+1}^{(1)}$ for $1\leq j\leq n-1$.
 \end{prop}

\begin{proof}
Let $\gb_n:=\set{c_j : 1\leq j\leq n}$ be a constrained optimal set of $n$-points for $\mu$ with respect to the constraint $L$. Then, $c_j=(a_j, ma_j+c)$, where $a_j\in \D R$ for $1\leq j\leq n$. Without any loss of generality, we can assume that $a_1<a_2<\cdots<a_n$. Let the Voronoi region of $c_j$ contains the set $[q_j^{(1)}, q_j^{(2)}]$ from the support of $\mu$ for $1\leq j\leq n$, where $q_1^{(1)}=1$, $q_n^{(2)}=\infty$ and $q_j^{(2)}+1=q_{j+1}^{(1)}$ for $1\leq j\leq n-1$.
Then, we have
\begin{align}\label{Megha43}
Q_n^{(c)}=\sum_{j=1}^{n} \sum_{i=q_j^{(1)}}^{q_j^{(2)}}\frac 1{2^i} \rh(i, c_j).
\end{align}
Notice that for $1\leq j\leq n$, the distortion error contributed by $c_j$ on its own Voronoi region is given by
\begin{align} \label{Megha44}
 &\sum_{i=q_j^{(1)}}^{q_j^{(2)}}\frac 1{2^i} \rh\Big(i, (a_j, ma_j+c)\Big) \notag\\
 &=\sum_{i=q_j^{(1)}}^{q_j^{(2)}}\frac 1{2^i}\Big((i -a_j)^2+( ma_j+c)^2\Big)\notag\\
 &=\sum_{i=q_j^{(1)}}^{q_j^{(2)}}\frac 1{2^i} \Big(\Big((i - Av[q_j^{(1)}, q_j^{(2)}])+(Av[q_j^{(1)}, q_j^{(2)}]-a_j)\Big)^2 +(ma_j+c)^2\Big)\notag\\
&=\sum_{i=q_j^{(1)}}^{q_j^{(2)}}\frac 1{2^i} \Big( (i - Av[q_j^{(1)}, q_j^{(2)}] )^2+2(i - Av[q_j^{(1)}, q_j^{(2)}] )( Av[q_j^{(1)}, q_j^{(2)}]-a_j)\notag\\
&\qquad \qquad + ( Av[q_j^{(1)}, q_j^{(2)}]-a_j)^2 +(ma_j+c)^2\Big)\notag\\
&=\sum_{i=q_j^{(1)}}^{q_j^{(2)}}\frac 1{2^i} \Big((i -Av[q_j^{(1)}, q_j^{(2)}])^2+( Av[q_j^{(1)}, q_j^{(2)}]-a_j)^2 +(ma_j+c)^2\Big).
\end{align}
Notice that the above distortion error is minimum if $( Av[q_j^{(1)}, q_j^{(2)}]-a_j)^2 +(ma_j+c)^2$ is minimum, i.e., if
\[a_j=\frac{Av[q_j^{(1)}, q_j^{(2)}]-cm}{1+m^2},\]
and then
\begin{equation} \label{Megha45} c_j=\Big(\frac{Av[q_j^{(1)}, q_j^{(2)}]-cm}{1+m^2}, \frac{ mAv[q_j^{(1)}, q_j^{(2)}]+c}{1+m^2}\Big)=U(Av[q_j^{(1)}, q_j^{(2)}]), \end{equation}
where $1\leq j\leq n$. Putting $a_j=\frac{Av[q_j^{(1)}, q_j^{(2)}]-cm}{1+m^2}$ in \eqref{Megha44}, we have
\[\sum_{i=q_j^{(1)}}^{q_j^{(2)}}\frac 1{2^i} \rh\Big(i, (a_j, ma_j+c)\Big)=\sum_{i=q_j^{(1)}}^{q_j^{(2)}}\frac 1{2^i} \Big((i -Av[q_j^{(1)}, q_j^{(2)}])^2+\frac{(m Av[q_j^{(1)}, q_j^{(2)}] +c)^2}{1+m^2}\Big).\]
Hence, by \eqref{Megha43}, we have
\begin{align*}
Q_n^{(c)}=\sum_{j=1}^{n}\sum_{i=q_j^{(1)}}^{q_j^{(2)}}\frac 1{2^i} (i -Av[q_j^{(1)}, q_j^{(2)}])^2+\frac 1 {1+m^2}\sum_{j=1}^{n}(m Av[q_j^{(1)}, q_j^{(2)}] +c)^2(2^{1-q_j^{(1)}}-2^{-q_j^{(2)}}).
\end{align*}
Let us now consider the set $\ga_n$ given by
\begin{equation} \label{Megha46}\ga_n=\Big\{Av[q_j^{(1)}, q_j^{(2)}] : 1\leq j\leq n\Big\}.\end{equation}
Notice that the unconstrained distortion error for $\mu$ with respect to the set $\ga_n$, denoted by $Q(\mu; \ga_n)$, is defined as
\[Q(\mu; \ga_n)=\sum_{j=1}^{n} \sum_{i=q_j^{(1)}}^{q_j^{(2)}}\frac 1{2^i}(i -Av[q_j^{(1)}, q_j^{(2)}])^2.\]
Thus, we see that
\begin{equation} \label{Megha47} Q_n^{(c)}=Q(\mu; \ga_n)+\frac 1 {1+m^2}\sum_{j=1}^{n}(m Av[q_j^{(1)}, q_j^{(2)}] +c)^2(2^{1-q_j^{(1)}}-2^{-q_j^{(2)}}).
\end{equation}
From the above expression, we see that $Q_n^{(c)}$ will be minimum iff the distortion error $Q(\mu; \ga_n)$ is minimum with respect to all possible subsets of $\te{supp}(\mu)$ with cardinality $n$, i.e., when $Q(\mu; \ga_n)$ becomes $Q_n^{(u)}$. Thus, we can say that the set $\gb_n$  will be a constrained optimal set of $n$-points for $\mu$ iff the set $\ga_n$ forms an unconstrained optimal set of $n$-means for $\mu$. By \eqref{Megha45} and \eqref{Megha46}, we have $\gb_n=U(\ga_n)$, and then as a result of \eqref{Megha47},
the proof of the proposition is complete.
\end{proof}

 The following theorems give the main results in this section.

 \begin{theorem} \label{Theo2}
Let $\gb_n$ be a constrained optimal set of $n$-points for any $n\in \D N$ with $n$th constrained quantization error $Q_n^{(c)}$. Then,
 \begin{enumerate}
 \item $\gb_1=\Big\{(\frac{2-c m}{m^2+1},\frac{c+2 m}{m^2+1})\Big\} \te{ with } Q_1^{(c)}=\frac{c^2+4 c m+6 m^2+2}{m^2+1}$.
 \item $\gb_2=\Big\{(\frac{4-3 c m}{3 m^2+3},\frac{c+\frac{4 m}{3}}{m^2+1}),  (\frac{4-c m}{m^2+1},\frac{c+4 m}{m^2+1})\Big\} \te{ with } Q_2^{(c)}=\frac{3 c^2+12 c m+18 m^2+2}{3 m^2+3}$.
     \item $\gb_3=\Big\{(\frac{1-c m}{m^2+1},\frac{c+m}{m^2+1}),(\frac{7-3 c m}{3 m^2+3},\frac{c+\frac{7 m}{3}}{m^2+1}),(\frac{5-c m}{m^2+1},\frac{c+5 m}{m^2+1})\Big\} \\
         \te{ or }  \Big\{ (\frac{4-3 c m}{3 m^2+3},\frac{c+\frac{4 m}{3}}{m^2+1}),(\frac{10-3 c m}{3 m^2+3},\frac{c+\frac{10 m}{3}}{m^2+1}),(\frac{6-c m}{m^2+1},\frac{c+6 m}{m^2+1})\Big\} \te{ with } Q_3^{(c)}=\frac{3 c^2+12 c m+18 m^2+1}{3 m^2+3}$.

\item For $n\geq 4$, \\
$\gb_n=\Big\{(\frac{j-c m}{m^2+1},\frac{c+j m}{m^2+1}) : 1\leq j\leq n-2\Big\}\UU \Big\{(-\frac{3 c m-3 n+2}{3 m^2+3},\frac{3 c+3 m n-2 m}{3 m^2+3})\Big\}
\UU\Big\{(\frac{-c m+n+2}{m^2+1},\frac{c+m (n+2)}{m^2+1})\Big\} \\
\quad \te{ or }\\
\gb_n=\Big\{(\frac{j-c m}{m^2+1},\frac{c+j m}{m^2+1}) : 1\leq j\leq n-3\Big\}\UU \Big\{(-\frac{3 c m-3 n+5}{3 m^2+3},\frac{3 c+3 m n-5 m}{3 m^2+3})\Big\}\UU\Big\{(\frac{-3 c m+3 n+1}{3 m^2+3},\frac{3 c+3 m n+m}{3 m^2+3})\Big\}\\
\hspace{ 1 in} \UU \Big\{(\frac{-c m+n+3}{m^2+1},\frac{c+m (n+3)}{m^2+1})\Big\} \\
\te{ with } Q_n^{(c)}=\frac{3 c^2+12 c m+18 m^2+2^{3-n}}{3 m^2+3}$.
\end{enumerate}
 \end{theorem}
 \begin{proof}
Let $\ga_n$ be an unconstrained optimal set of $n$-means for any $n\in \D N$ with $n$th unconstrained quantization error $Q_n^{(u)}$. Then, by \cite[Section~4]{CHMR}, it is known that
\begin{enumerate}
 \item $\ga_1=\Big\{Av[1, \infty]\Big\} \te{ with } Q_1^{(u)}=Er[1, \infty]=2$.
 \item $\ga_2=\Big\{Av[1, 2], Av[3, \infty]\Big\} \te{ with } Q_2^{(u)}=Er[1, 2]+Er[3, \infty]=\frac 23$.
     \item $\ga_3=\Big\{1, Av[2,3], Av[4, \infty]\Big\}, \te{ or } \Big\{Av[1,2], Av[3, 4], Av[5, \infty]\Big\}  \te{ with } Q_3^{(u)}=\frac 13.$

\item For $n\geq 4$,
$\ga_n=\Big\{1, 2, 3, \cdots, n-2, Av[n-1,n], Av[n+1, \infty]\Big\}, \\
\te{ or } \Big\{1, 2, 3, \cdots, n-3, Av[n-2,n-1], Av[n, n+1],  Av[n+2, \infty]\Big\}\\
\te{ with } Q_n^{(u)}=\frac{2^{3-n}}{3}$.
\end{enumerate}
Then, by Proposition~\ref{propim1}, the results of the theorem can be obtained. However, for the readers' convenience, we give a detailed proof for $n \geq 4$.
Choose $\ga_n=\Big\{1, 2, 3, \cdots, n-2, Av[n-1,n], Av[n+1, \infty]\Big\}$. Then,
\begin{align*}
\gb_n&=U(\ga_n)=\set{U(x) : x\in \ga_n}\\
&=\Big\{(\frac{j-c m}{m^2+1},\frac{c+j m}{m^2+1}) : 1\leq j\leq n-2\Big\}\UU \Big\{\big(-\frac{3 c m-3 n+2}{3 m^2+3},\frac{3 c+3 m n-2 m}{3 m^2+3}\big)\Big\}\\
&\quad \qquad \qquad  \UU\Big\{(\frac{-c m+n+2}{m^2+1},\frac{c+m (n+2)}{m^2+1})\Big\}.
 \end{align*}
By \eqref{Megha42}, the $n$th constrained quantization is given by
\begin{align*}
 Q_n^{(c)}&=Q_n^{(u)}+\frac 1 {1+m^2}\Big(\sum_{j=1}^{n-2} (mj+c)^2 \frac 1{2^j} +(m Av[n-1, n]+c)^2(2^{1-(n-1)}-2^{-n})\\
 &\qquad \qquad +(m Av[n+1, \infty]+c)^2 2^{1-(n+1)}\Big)\\
 &=\frac{3 c^2+12 c m+18 m^2+2^{3-n}}{3 m^2+3}.
\end{align*}
If we choose $\ga_n=\Big\{1, 2, 3, \cdots, n-3, Av[n-2,n-1], Av[n, n+1],  Av[n+2, \infty]\Big\}$, similarly we can show that
$\gb_n=\Big\{(\frac{j-c m}{m^2+1},\frac{c+j m}{m^2+1}) : 1\leq j\leq n-3\Big\}\UU \Big\{(-\frac{3 c m-3 n+5}{3 m^2+3},\frac{3 c+3 m n-5 m}{3 m^2+3})\Big\}\UU\Big\{(\frac{-3 c m+3 n+1}{3 m^2+3},\frac{3 c+3 m n+m}{3 m^2+3})\Big\}\\
\hspace{ 1 in} \UU \Big\{(\frac{-c m+n+3}{m^2+1},\frac{c+m (n+3)}{m^2+1})\Big\}$
with
\begin{align*}
 Q_n^{(c)}&=Q_n^{(u)}+\frac 1 {1+m^2}\Big(\sum_{j=1}^{n-3} (mj+c)^2 \frac 1{2^j} +(m Av[n-2, n-1]+c)^2(2^{1-(n-2)}-2^{-(n-1)})\\
 &\qquad \qquad+(m Av[n, n+1]+c)^2(2^{1-n}-2^{-(n+1)}) +(m Av[n+2, \infty]+c)^2 2^{1-(n+2)}\Big)\\
 &=\frac{3 c^2+12 c m+18 m^2+2^{3-n}}{3 m^2+3}.
\end{align*}
Thus, the proof of the theorem is complete.
 \end{proof}

 \begin{theorem} \label{Theo3}
 Let $\mu$ be the discrete probability measure defined by $\mu= \sum_{n=1}^\infty \frac 1 {2^n} \, \delta_n$. Then, with respect to the constraint $L$ given by $y=mx+c$, where $m, c\in \D R$, the constrained quantization dimension is zero.
 \end{theorem}

 \begin{proof}
 Let $Q_n^{(c)}(\mu)$ be the $n$th constrained quantization error with respect to the constraint $L$. Then,
\[Q^{(c)}_{\infty} (\mu)=\lim_{n\to \infty} Q_n^{(c)}(\mu)=\frac{c^2+4 c m+6 m^2}{m^2+1}.\]
Hence,
\[\lim_{n\to \infty} \frac{2\log n}{-\log (Q^{(c)}_n(\mu)-Q^{(c)}_\infty(\mu))}=0\]
yielding the fact that the constrained quantization dimension is zero.
 \end{proof}

 \begin{remark} \label{remark467}
Theorem~\ref{Theo3} implies that for the discrete distribution $\mu$ defined by
\[
\mu = \sum_{n=1}^\infty \frac{1}{2^n} \, \delta_n,
\]
with respect to the constraint $L$ given by the line $y = mx + c$, where $m, c \in \mathbb{R}$, the constrained quantization dimension is zero. In particular, when $m = c = 0$, the constraint $L$ reduces to the real line $\mathbb{R}$, and hence the unconstrained quantization dimension for the infinite distribution $\mu$ is also zero.
By~\cite{PR1}, it is known that for a Borel probability measure, the constrained quantization dimension and the constrained quantization coefficient depend on the underlying constraint. Therefore, it is natural to investigate whether there exists a constraint $L$ for which the constrained quantization dimension of the discrete distribution $\mu$ differs from the unconstrained quantization dimension. For a certain constraint if the constrained quantization dimension $D(\mu)$ exists as a finite positive number, then it will be worthwhile to investigate whether $D(\mu)$-dimensional constrained quantization coefficient exists as a finite positive number.
\end{remark}

\section{Conclusion and Future Work} \label{sec5}

This paper provides a detailed investigation into constrained quantization for both finite and infinite discrete probability distributions supported on subsets of the real line. Under two primary geometric constraints-a semicircular arc and the union of two sides of an equilateral triangle, we determined the constrained optimal sets of $n$-points and computed the corresponding constrained quantization errors for various distributions.

For finite discrete distributions, we analyzed both uniform and nonuniform distributions supported on the finite set $\{-3, -2, \ldots, 3\}$. For each case and constraint, we derived explicit expressions for the constrained quantization errors and identified all possible constrained optimal sets for small values of $n$, often accompanied by multiple canonical vector configurations.

For infinite discrete distributions, two distinct cases were considered: one supported on the set $\left\{\frac{1}{n} : n \in \mathbb{N} \right\}$ for which we computed constrained optimal sets and quantization errors for moderate values of $n$, and the other on the set $\mathbb{N}$ for which we developed a general theoretical framework for constrained quantization  with respect to a linear constraint $L: y = mx + c$. The framework established a direct connection between constrained and unconstrained quantizers via bijective geometric transformations, and derived an exact formula relating constrained and unconstrained quantization errors.

\medskip

\noindent\textbf{Future Work.} Several natural extensions of this work merit further investigation:

\begin{itemize}
    \item Extend the theory to discrete distributions with more general or multidimensional supports, including those exhibiting polynomial or exponential decay.
    \item Investigate constrained quantization under more complex geometric constraints such as curves, manifolds, or fractal boundaries in higher-dimensional Euclidean spaces.
    \item Analyze the limiting behavior of constrained quantization errors and rigorously determine the existence and values of constrained quantization dimensions and coefficients.
    \item Explore conditions under which the constrained quantization dimension exceeds, equals, or falls below its unconstrained counterpart, particularly for singular discrete measures.
    \item Develop efficient computational algorithms for determining constrained optimal sets in high-complexity or large-$n$ regimes, and explore their numerical stability and convergence.
\end{itemize}

\noindent The results in this paper reveal that geometric constraints significantly influence both the structure and existence of optimal quantizers, opening a broad avenue for continued research in constrained quantization theory.

\section*{Declaration}
							
							\noindent
\textbf{Authors' contributions:} Each author contributed equally to this manuscript. All authors have read and agreed to the published version of the manuscript.\\
\\
\noindent \tbf{Funding:}  This research received no external funding.\\
							
							\noindent
							\textbf{Data availability:} No data were used to support this study.\\
							\\
							\noindent
\textbf{Conflicts of interest.} The authors declare no conflict of interest.\\

\end{document}